\documentclass[11pt]{amsart}

\usepackage{amssymb}

\usepackage{url}
\usepackage[all]{xy}

\usepackage[colorlinks=true,linkcolor=blue,urlcolor=blue,citecolor=blue,pagebackref]{hyperref} 

\usepackage[usenames, dvipsnames, svgnames]{xcolor}

\usepackage{enumerate} 

\usepackage{helvet} 
\usepackage{courier} 
\usepackage{eucal} 
\usepackage{mathrsfs} 

\reversemarginpar 
\normalmarginpar 

\let\oldmarginpar\marginpar 
\renewcommand\marginpar[1]{\-\oldmarginpar{\raggedright\small\sf #1}}

\title[Finite group scheme actions and incompressibility of
covers]{Finite groups scheme actions and incompressibility of Galois
  covers: beyond the ordinary case.}

\author{Najmuddin Fakhruddin}
\author{Rijul Saini} 
\address{School of Mathematics, Tata Institute of Fundamental Research, 
Homi Bhabha Road, Mumbai 400005, India}

\email{naf@math.tifr.res.in}
\email{rijul@math.tifr.res.in}

\newcommand{\nc}{\newcommand}

\nc{\rnc}{\renewcommand}

\nc{\bs}{\backslash}
\nc{\te}{\otimes}
\nc{\lf}{\lfloor} 
\nc{\rf}{\rfloor}
\nc{\lc}{\lceil}  
\nc{\rc}{\rceil}
\nc{\lr}{\longrightarrow}
\nc{\sr}{\stackrel}
\nc{\dar}{\dashrightarrow}
\nc{\thra}{\twoheadrightarrow}

\nc{\la}{\langle}
\nc{\ra}{\rangle} 

\nc{\ms}{\mathscr}
\nc{\mc}{\mathcal}
\nc{\mb}{\mathbb}
\nc{\mf}{\mathbf}
\nc{\mr}{\mathrm}
\nc{\mg}{\mathfrak}
\nc{\msf}{\mathsf}

\nc{\bP}{\mathbb{P}}
\rnc{\P}{\mathbb{P}}
\nc{\Q}{\mathbb{Q}}
\nc{\Z}{\mathbb{Z}}
\nc{\C}{\mathbb{C}}
\nc{\R}{\mathbb{R}}
\nc{\A}{\mathbb{A}}
\nc{\V}{\mathbb{V}}
\nc{\W}{\mathbb{W}}
\nc{\N}{\mathbb{N}}
\nc{\D}{\mathbb{D}}
\nc{\G}{\mathbb{G}}
\nc{\F}{\mathbb{F}}

\nc{\X}{\mathbb{X}}
\nc{\Y}{\mb{Y}}
\nc{\qb}{\overline{\mathbb{Q}}}

\nc{\del}{\partial}

\nc{\wt}{\widetilde}
\nc{\wh}{\widehat}
\nc{\ov}{\overline}
\nc{\un}{\underline}

\nc{\aff}{{\A}^1}
\nc{\naive}{\!\sim_n}

\nc{\Spec}{\mr{Spec}}

\nc{\ord}{\mr{ord}}

\nc{\omx}{\omega_X}

\nc{\ep}{\epsilon}
\nc{\ve}{\varepsilon}
\nc{\vt}{\vartheta}

\rnc{\l}{\lambda}
\rnc{\k}{\kappa}
\nc{\tk}{\tilde{\kappa}}
\rnc{\d}{\delta}

\nc{\can}{\mr{can}}

\nc{\ovl}{\ov{\lambda}}
\nc{\vl}{\mb{V}_{\ovl}}
\nc{\dl}{\mb{D}_{\ovl}}
\nc{\mnb}{\ov{\mr{M}}_{0,n}}
\nc{\mn}{\mr{M}_{0,n}}
\nc{\mel}{\ov{\mr{M}}_{1,1}}
\nc{\mfb}{\ov{\mr{M}}_{0,4}}
\nc{\mof}{\mr{M}_{0,4}}
\nc{\mgnb}{\ov{\mr{M}}_{g,n}}
\nc{\mgn}{\ov{\mr{M}}_{g,n}}
\nc{\omc}{\ov{\mr{M}}}

\rnc{\sl}{\shoveleft}

\nc{\res}{\operatorname{Res}}
\nc{\pic}{\operatorname{Pic}}
\nc{\spec}{\operatorname{Spec}}
\nc{\im}{\operatorname{Im}}
\nc{\gal}{\operatorname{Gal}}
\nc{\fr}{\operatorname{Fr}}
\nc{\ed}{\operatorname{ed}}
\nc{\rank}{\operatorname{rank}}
\nc{\h}{\operatorname{H}}
\nc{\ch}{\operatorname{char}}
\nc{\sw}{\operatorname{sw}}
\nc{\rsw}{\operatorname{rsw}}
\nc{\supp}{\operatorname{supp}}

\nc{\stab}{\operatorname{Stab}}

\nc{\Mor}{\operatorname{Mor}}
\nc{\Per}{\operatorname{Per}}
\nc{\prep}{\operatorname{Prep}}
\nc{\End}{\operatorname{End}}
\nc{\Orb}{\operatorname{Orb}}

\nc{\Aut}{\mf{{Aut}}}
\nc{\UAut}{\mf{UAut}}
\nc{\Hom}{\mf{{Hom}}}
\nc{\Ext}{\mf{Ext}}

\rnc{\End}{\mf{End}}
\rnc{\R}{R}
\nc{\T}{\mc{T}}

\rnc{\hom}{\operatorname{Hom}}
\nc{\en}{\operatorname{End}}
\nc{\aut}{\operatorname{Aut}}

\rnc{\G}{G}

\nc{\f}{f}
\nc{\ff}{\mg{f}}
\nc{\bt}{\bar{\tau}}

\newtheorem{thm}{Theorem}[section]
\newtheorem{prop}[thm]{Proposition}
\newtheorem{conj}[thm]{Conjecture}
\newtheorem{cor}[thm]{Corollary}
\newtheorem{lem}[thm]{Lemma}

\theoremstyle{definition}

\newtheorem{ques}[thm]{Question}
\theoremstyle{remark}
\newtheorem{rem}[thm]{Remark}

\newtheorem{ex}[thm]{Example}

\numberwithin{equation}{section}

\begin{document}

\begin{abstract}
  Inspired by recent work of Farb, Kisin and Wolfson \cite{FKW}, we
  develop a method for using actions of finite group schemes over a
  mixed characteristic dvr $\R$ to get lower bounds for the essential
  dimension of a cover of a variety over $K = \mr{Frac}(\R)$. We then
  apply this to prove $p$-incompressibility for congruence covers of a
  class of unitary Shimura varieties for primes $p$ at which the
  reduction of the Shimura variety (at any prime of the reflex field
  over $p$) does not have any ordinary points. We also make some
  progress towards a conjecture of Brosnan on the
  $p$-incompressibility of the multiplication by $p$ map of an abelian
  variety.
\end{abstract}

\maketitle

\section{Introduction} \label{s:intro}

Let $f:X \to Y$ be a generically etale morphism of varieties over a
field $K$ with $Y$ integral. We say that $f$ is \emph{incompressible}
if there is no commutative Cartesian diagram
\[
  \xymatrix{ X \ar@{.>}[r] \ar[d]^f  &X' \ar[d]^{f'} \\
    Y \ar@{.>}[r] & Y'
  }
\]
with $Y'$ integral, $\dim(Y') < \dim(Y)$, and the horizontal arrows
are dominant rational maps. If $p$ is a prime number, we say that $f$
is $p$-\emph{incompressible} if the morphism $X \times_Y Z \to Z$ is
incompressible, where $Z$ is integral and $Z \to Y$ is a dominant
generically finite morphism of degree prime to $p$. For more on this
definition (and the more general notion of \emph{essential dimension}
of $f$) the reader may consult \cite{buhler-reichstein}, \cite{FKW} or
\cite{merk-survey}.

In the article \cite{FKW}, Farb, Kisin and Wolfson introduced a new
method for proving $p$-incompressibility over fields of characteristic
zero by using an elegant mixed characteristic method and applied it to
prove the $p$-incompressibility of a large class of congruence covers
of Shimura varieties of Hodge type under an \emph{ordinarity}
assumption on the reduction of the Shimura variety at a prime of the
reflex field lying over $p$.  Their result was reproved in many cases,
and extended to certain Shimura varieties not of Hodge type, using the
``fixed-point method'' by Brosnan and the first-named author in
\cite{bf-fixed}. Moreover, the latter authors formulated a general
conjecture \cite[Conjecture 1]{bf-fixed} on the essential dimension of
covers associated to arbitrary variations of Hodge structure which, in
the context of Shimura varieties, predicts that the ordinarity
assumption of \cite{FKW} is unnecessary.  The main goal of this paper
is to develop methods which allow us to prove that this is indeed true
for a certain class of unitary Shimura varieties (Theorem
\ref{t:main}) and using related methods we also make progress towards
a conjecture of P.~Brosnan (Conjecture \ref{c:abelian}) on the
$p$-incompressibility of the multiplication by $p$ map $[p]:A \to A$,
where $A$ is an abelian variety over a field of characteristic zero
(Theorem \ref{t:ab}).

A special case of our theorem on unitary Shimura varieties is the
following result, which we state somewhat informally here. For the
precise (and more general) statement, see Theorem \ref{t:main}.
\begin{thm} \label{t:quadratic}%
  Let $F$ be an imaginary quadratic extension of $\Q$ and let $S$ be a
  PEL moduli space of abelian varieties of odd dimension
  $d = 2\d + 1\geq 3$ with endomorphism ring containing $O_F$. If the type of
  the $O_F$-action is $(\d, \d+1)$, then for all primes of good
  reduction of $S$ the principal level $p$ congruence cover
  $S(p) \to S$ is $p$-incompressible.
\end{thm}
If one specializes \cite[Theorem 4.3.6]{FKW} to this case, one obtains
$p$-incompressibility only for primes $p$ which split in
$O_F$. Theorem \ref{t:main} applies to PEL Shimura varieties
corresponding to arbitrary CM fields $F$ (with a similar restriction
on the ``type'') and though in this setting our results are strictly
stronger that those of \cite{FKW}, they require that $p$ split
completely in the maximal totally real subfield $F_0$ of $F$.

A consequence of our main theorem (Theorem \ref{t:ab}) on abelian
varieties is the following (Corollary \ref{c:three}):
\begin{thm} \label{t:three}%
  For any abelian variety $A$ of dimension $d \leq 3$ over a field of
  characteristic zero there exists a set $\mg{P}(A)$ of rational
  primes of positive density such that $[p]:A \to A$ is
  $p$-incompressible for $p \in \mg{P}(A)$.
\end{thm}

Our methods, which are inspired by and extend those of Farb, Kisin and
Wolfson, depend on degenerating the covers to fields of characteristic
$p$.  The two key ingredients used in \cite{FKW} are the Serre--Tate
theorem relating deformations of abelian varieties and their
$p$-divisible groups, and Kummer theory, which gives a very concrete
description of $\mu_p$-torsors (over any local base). At an ordinary
point, the Serre--Tate theorem gives a very explicit (formal)
description of the degeneration of the congruence cover when one
specializes to characteristic $p$ and Kummer theory allows the authors
of \cite{FKW} to prove an extension result (\cite[Lemma 3.1.5]{FKW})
for $(\mu_p)^n$-torsors which plays an essential role in their proofs
of $p$-incompressibility.

For torsors under more general finite (possibly noncommutative) group
schemes there is no good analogue of Kummer theory, so we change our
point of view and consider actions of finite group schemes rather than
torsors. In this context, our replacement for \cite[Lemma 3.1.5]{FKW}
is Lemma \ref{l:descent}, an extension result for actions of general
finite flat group schemes. Not working with torsors has drawbacks
though, since general actions of finite group schemes, even if they
are faithful (which is not guaranteed), can have complicated
stabilizers, and this is the main reason why we cannot prove
incompressibilty in the full expected generality (in the settings we
consider). In \S\ref{s:gsa} we consider two classes of group schemes
for which Lemma \ref{l:descent} can be usefully applied to give lower
bounds for essential dimension of covers; the first class is applied
to prove our results on Shimura varieties and the second on abelian
varieties. The main idea in both cases is the simple fact that we can
bound from below the dimension of a smooth variety on which a (finite)
group scheme acts using knowledge of the dimension of the Lie algebra
of the group scheme and that of the stabilizer of a general point.

The group schemes used for the application to Shimura varieties are
noncommutative (in the non-ordinary case, which is our main interest)
and \S\ref{s:mu} is devoted to an analysis of the structure of these
group schemes in the setting of PEL Shimura varieties of type A. Here
we depend heavily on some results of Moonen from
\cite{moonen-serre-tate}. In \S\ref{s:shim} we construct certain
integral models of the $p$-congruence cover of such Shimura varieties
(for $p$ an unramified prime) by a naive method, i.e., by
normalisation of the Kottwitz integral model \cite{kottwitz-points} in
the function field of the $p$-congruence cover. We show that these
integral models are generically smooth along the special fibre and
admit an action of a suitable finite flat group scheme using a result
of Cariani--Scholze from \cite{cs-compact-unitary} and the results of
\S\ref{s:mu}. In \S\ref{s:ed} we apply Proposition \ref{p:incrit} (our
incompressibilty criterion for certain noncommutative group scheme
action) to this integral model to prove our main result on Shimura
varieties, Theorem \ref{t:main}. Finally, in \S\ref{s:abelian} we
prove our results on abelian varieties. Here we use commutative group
schemes which are ``almost ordinary'' via Proposition \ref{p:rank}, a
criterion giving a lower bound on the essential dimension for actions
of such group schemes.

\subsection{Notation}

We usually denote by $K$ a complete discretely valued field of
characteristic $0$, by $\R$ its ring of integers and by $k$ its
residue field which will always be perfect (and usually, but not
always, of characteristic $p>0$). We denote $\spec(\R)$ by $\T$ and
usually (but not always) denote schemes over $\T$ by using
calligraphic fonts, e.g., $\mc{X}$, $\mc{G}$, \ldots, and the generic
fibre of such a scheme by the corresponding letter in ordinary font,
i.e., $X =\mc{X}_K$, $G = \mc{G}_K$, \ldots, and we will say that
$\mc{X}$ (resp. $\mc{G}$) is a model of $X$ (resp. $G$), \ldots. We
will implicitly assume (except when explicitly mentioned otherwise)
that the scheme $\mc{X}$ over $\T$ is flat and the closed fibre, which
will usually be denoted by using a subscript $0$, e.g., $\mc{X}_0$,
$\mc{G}_0$, \ldots, is non-empty.  However, we will often use notation
such as $\mc{Z}_0$ to denote a subscheme of $\mc{X}_0$ as above
without assuming that it is the closed fibre of a flat subscheme
$\mc{Z}$ of $\mc{X}$.

For any finite connected (not necessarily flat) non-empty scheme
$\mc{X}$ over $\T$ we will denote by $T(\mc{X})$ its reduced
tangent space, i.e., the tangent space of $\mc{X}_0$ (at its unique
point). We set $t(\mc{X}) = \dim(T(\mc{X}))$.

For a finite commutative group scheme (resp.~$p$-divisible group) over
any base, we use $D$ as a superscript (e.g., $\mc{G}^D$) to denote its
Cartier (resp.~Serre) dual. We also recall that a finite flat group
scheme (over any base) is \emph{multiplicative} if its Cartier dual is
(finite) etale.

\subsection{Acknowledgements}

We thank Patrick Brosnan, Arvind Nair, Madhav Nori and Vincent Pilloni
for useful conversations and/or correspondence and Zinovy Reichstein
for his comments on a preliminary version of this paper.  N.F. was
supported by the DAE, Government of India, under P.I.No.~RTI4001.

\section{Finite group scheme actions and essential
  dimension} \label{s:gsa}%

\subsection{Descending group scheme actions}

The lemma below is the key property of multiplicative group schemes
that (from our point of view) makes them very useful for proving
incompressibility of covers. Although the goal of this paper is to
prove incompressibility using more general group schemes,
multiplicative group schemes (as subgroup schemes of more general
group schemes) still play a key role in all our results.

\begin{lem} \label{l:mult}%
  Let $\mc{G}$ be a finite flat connected multiplicative group scheme
  over $\T$. If $\mc{G}$ acts on a scheme $\mc{X}$ which is smooth and
  of finite type over $\T$ (with $\mc{X}_0 \neq \emptyset$) so that
  $G$ acts generically freely on $X$, then $\mc{G}$ acts freely on an
  open subscheme $\mc{X}' \subset \mc{X}$ (with
  $\mc{X}'_0 \neq \emptyset$).
\end{lem}

\begin{proof}
  It suffices to show that $\mc{G}_0$ acts generically freely on
  $\mc{X}_0$, so after a base change we may assume that $k$ is
  algebraically closed which implies that $\mc{G}^D$ is constant. For
  the action of a finite group scheme freeness of the action is
  equivalent to the stabilizers of all geometric points being trivial,
  so if generic freeness does not hold then since $\mc{X}_0$ is smooth
  and $\mc{G}_0$ only has finitely many subgroup schemes it follows
  that there exists a non-empty open set of $\mc{X}_0$ such that all
  geometric points of this subset have the same non-trivial
  stabilizer. We choose any $\mc{G}$-invariant affine open
  $\mc{X}' \subset \mc{X}$ (with $\mc{X}'_0 \neq \emptyset$) such that
  all geometric points of $\mc{X}'_0$ have the same non-trivial
  stabilizer in $\mc{G}_0$.  Since $\mc{G}^D$ is constant, one then
  easily reduces to the case that $\mc{G} = \mu_p$ and the action of
  $\mc{G}_0$ on $\mc{X}'_0$ is trivial; here we use smoothness of
  $\mc{X}_0'$ to conclude that if $\mu_p$ stabilizes each geometric
  point then the action must be trivial.

  Writing $\mc{X}'_0 = \spec(B)$, the $\mu_p$ action translates into
  the structure of a $\Z/p\Z$ grading on the $\R$-algebra $B$, i.e.,
  $B = \oplus_{i \in \Z/p\Z} B_i$ with
  $B_i \cdot B_j \subset B_{i+j}$. Since this grading is compatible
  with base change, the assumption that the action on $\mc{X}'_0$ is
  trivial implies that $B_i \otimes_{\R} k = \{0\}$ for $i \neq 0$.
  This implies that the $G$-action on the image of $\mc{X}'(R)$ in
  $X(K)$ is trivial, but the assumptions on $\mc{X}$ and $\mc{T}$
  imply that this set is Zariski dense in $X$, thereby contradicting
  the generic freeness of the $G$-action on $X$.
\end{proof}

  \begin{lem} \label{l:descent}%
  Let $\mc{G}$ be a finite flat
  group scheme over $\T$ acting on a separated irreducible normal
  scheme $\mc{X}$, faithfully flat and of finite type over $\T$, with
  $\mc{G}_0(k)$ acting trivially on the set of irreducible components
  of $\mc{X}_0$. Assume further that $G = \mc{G}_K$ is a constant
  group scheme and there is a $G$-equivariant dominant morphism
  $f: X \to Y$, where $Y$ is a smooth quasi-projective variety over
  $K$. Then there exist
  \begin{enumerate}
  \item a non-empty open $G$-invariant subvariety $Y' \subset Y$ and a
    normal model $\mc{Y}'$ of $Y'$ over $\T$ (with
    $\mc{Y}'_0 \neq \emptyset$) such that the $G$-action on $Y'$
    extends to a $\mc{G}$-action on $\mc{Y}'$ and
  \item a non-empty open $\mc{G}$-invariant subscheme $\mc{X}' \subset
    \mc{X}$,
  \end{enumerate}
  such that $f|_{X'}$ extends to a $\mc{G}$-equivariant surjective
  flat morphism $\phi: \mc{X}' \to
  \mc{Y}'$. 
  If $\mc{X}$ is smooth over $\T$ then we may also take $\mc{Y}'$
  to be smooth over $\T$.
\end{lem}

\begin{proof}
  Let $Z' = Y/G$ and let $\mc{Z}$ be a normal projective (flat) scheme
  over $\T$ containing $Z'$ as an open subscheme. Since $\mc{X}$ is
  normal and $\mc{Z}$ is proper, the rational map
  $g: \mc{X} \dar \mc{Z}$, induced by $f$ and the quotient map
  $Y \to Z'$, is defined at the generic points of $\mc{X}_0$.  By
  replacing $\mc{Z}$ by a suitable noramlized blow up of $\mc{Z}$
  along a closed subscheme of $\mc{Z}_0$, we may assume that there is
  an open $\mc{G}$-invariant subscheme $\mc{X}'$ of $\mc{X}$ on which
  $g$ is defined and the morphism $\mc{X}'_0 \to \mc{Z}_0$ dominates
  at least one irreducible component of $\mc{Z}_0$.

  Let $\tilde{\mc{Y}}$ be the normalisation of $\mc{Z}$ in the
  function field of $Y$. The group $G_{\T}$ (the group $G$ viewed as a
  constant group scheme over $\T$) acts on $\tilde{\mc{Y}}$ and it
  contains $Y$ as a $G$-invariant open subset. The morphism $f$
  extends to a $G_{\T}$-equivariant morphism
  $\phi: \mc{X}' \to \tilde{\mc{Y}}$ (since $\mc{X}'$ is normal,
  $\tilde{\mc{Y}}$ is finite over $\mc{Z}$ and its function field is
  contained in that of $\mc{X}'$) and $\mc{X}_0'$ dominates an
  irreducible component of $\tilde{\mc{Y}}_0$. By openness of flatness
  and the condition on the $\mc{G}_0(k)$-action on the irreducible
  components of $\mc{X}_0'$, there exists an open $G_{\T}$-invariant
  subscheme $\mc{Y}'$ of $\tilde{\mc{Y}}$ (with
  $\mc{Y}_0' \neq \emptyset$) such that $\phi$ restricted to $\mc{X}'$
  is flat over $\mc{Y}'$. By shrinking $\mc{Y}'$ further and replacing
  $\mc{X}'$ by $\phi^{-1}(\mc{Y}')$ we may assume that $Y' \subset Y$
  and the restriction of $\phi$ from $ \mc{X}' \to \mc{Y}'$ is flat,
  surjective and $G_{\T}$-equivariant.

  Let $O(\mc{G})$ denote the affine algebra of $\mc{G}$. The $\mc{G}$
  action on $\mc{X}'$ then corresponds to a map of sheaves
  $O_{\mc{X}'} \to O(\mc{G}) \otimes_{\R} O_{\mc{X}'}$ satisfying the
  usual identities. The $G_{\T}$-action on $\mc{Y}'$ induces a
  rational map $O_{\mc{Y}'} \dar O(\mc{G}) \otimes_{\R} O_{\mc{Y}'}$
  of sheaves, i.e., a map defined after tensoring with $K$.  The
  action of $\mc{G}$ extends to $\mc{Y}'$ if this map is in fact
  defined without tensoring with $K$ (since the identities needed for
  a group action will automatically hold by flatness over $\T$.)
  
  Since $O(\mc{G})$ is a finite free module, by using an
  $\R$-basis and the $G_{\T}$ equivariance of the map
  $\mc{X}' \to \mc{Y}'$ we are reduced to showing that if
  $\sigma \in \Gamma(Y',O_{Y'})$ is such that $f^*(\sigma)$ extends to
  an element of $\Gamma(\mc{X}', O_{\mc{X}'})$ then $\sigma$ extends
  to an element of $\Gamma(\mc{Y}', O_{\mc{Y}'})$. This follows
  immediately from the surjectivity of $\mc{X}' \to \mc{Y}'$ and the
  normality of $\mc{Y}'$ (see, e.g., \cite[Lemma 2.1]{nf-loc}).

  For the smoothness statement, since $k$ is perfect it suffices to
  observe that if $\mc{X}_0'$ is reduced then so is $\mc{Y}_0'$ since
  the map $\phi$ is flat.

\end{proof}

\begin{rem} \label{r:descent}%
  The lemma above is the basis for all our results on
  incompressibility. Note that although $\mc{G}$ is arbitrary, in this
  generality even if $\mc{G}$ acts freely on $\mc{X}$ and $Y$ we
  cannot conclude that it acts freely on $\mc{Y}'$ (see Example
  \ref{e:it} below). This is the fundamental difficulty one encounters
  when trying to use this lemma, and the rest of this section is
  devoted to developing methods which will allow us to overcome this
  difficulty in certain special cases.
\end{rem}

\begin{ex} \label{e:it}%
  Let $E$ be an elliptic curve over $K$ with good supersingular
  reduction and let $\mc{E}_{/\T}$ be its Neron model. Let $\mc{E}[p]$
  be the $p$-torsion subscheme of $\mc{E}$ and assume that its generic
  fibre is a constant group scheme (so isomorphic to
  $(\Z/p\Z)^2$). Let $\mc{G}', \mc{G}'' \subset \mc{E}[p]$ be distinct
  finite flat subgroup schemes of order $p$ and let
  $\mc{G} = \mc{G}'\times_{\T} \mc{G}''$. Clearly $t(\mc{G}) = 2$
  (since $\mc{E}_0$ is supersingular) and $\mc{G}$ acts freely on
  $\mc{E} \times_{\T} \mc{E}$. However, the sum map
  $\mc{E}^2 \to \mc{E}$ is $\mc{G}$-equivariant and $G$ acts freely on
  $E$ but there is a non-trivial subgroup scheme
  $\mc{H}_0 \subset \mc{G}_0$ which acts trivially on $\mc{E}_0$.
\end{ex}

\subsection{Group scheme actions, stabilizers and essential dimension
  of covers}

\begin{lem} \label{l:stab}%
  Let $\mc{G}$ be a finite flat connected group scheme acting on a
  smooth scheme $\mc{X}$ over $\T$. Suppose that for some closed point
  $x \in \mc{X}_0(k)$, $t(\mc{G}_0/\stab(x)) \geq n$. Then
  $\dim(X) \geq n$.
\end{lem}

\begin{proof}
  This is clear since the $\mc{G}$ action induces a closed embedding
  of $\mc{G}_0/\stab(x_0)$ into $\mc{X}_0$.  The hypothesis on the
  dimension and the smoothness of $\mc{X}_0$ implies that
  $\dim(\mc{X}_0) \geq n$ and so (by flatness) $\dim(X) \geq n$ as
  well.
\end{proof}

The proposition below allows us to use non-free $\mc{G}$-actions to
find lower bounds on the essential dimension of covers, but only for
very special $\mc{G}$. Note that for any non-trivial application
$\mc{G}$ has to be noncommutative.
  
\begin{prop} \label{p:incrit}%
  Let $G$ be a finite group acting faithfully on a smooth
  quasiprojective variety $X$ over $K$. Let $\mc{G}$ be a flat
  connected model of $G$ over $\T$ acting on a smooth model $\mc{X}$
  of $X$ over $\T$ (extending the $G$-action on $X$). Let $\mc{H}$ be
  a finite flat subgroup scheme of $\mc{G}$ which is of multiplicative
  type. Suppose that for all subgroup schemes $\mc{K}_0$ of $\mc{G}_0$
  which intersects $\mc{H}_0$ trivially, we have
  $t(\mc{G}_0/\mc{K}_0) \geq e$. Then
  \begin{enumerate}
  \item $\ed(X) \geq e$.
  \item If $\mc{G}$ acts freely on $\mc{X}$ then $\ed(X;p) \geq e$.
  \end{enumerate}
\end{prop}

\begin{proof}
  We first reduce (2) to (1), so let us assume that the $\mc{G}$
  action is free.  Then let $\mc{Z}$ be the quotient of $\mc{X}$ by
  $\mc{G}$; since $\mc{G}$ is finite and $\mc{X}$ is quasiprojective
  this exists (by the theory of Hilbert schemes) and the freeness of
  the $\mc{G}$-action implies that it is also smooth over $\T$; this
  can be checked fibre-wise and by \cite[\S12, Theorem 1]{av} the
  quotient map is flat which implies that if $\mc{X}$ is smooth then
  so is the quotient.
   Let $L/K(\mc{Z})$ be a finite extension of degree prime to $p$ and
  let $\mc{Z}'$ be the normalisation of $\mc{Z}$ in $L$. Using
  Abhyankar's lemma as in the proof of \cite[Theorem 3.2.6]{FKW},
   it follows that after replacing $K$ by a finite
  extension we may assume that the map $\mc{Z}' \to \mc{Z}$ is etale
  at at least one generic point of $\mc{Z}'_0$. By shrinking $\mc{Z}'$
  we may assume that the map $\mc{Z}' \to \mc{Z}$ is etale and, in
  particular, $\mc{Z}'$ is smooth. Let
  $\tilde{\mc{X}} = \mc{X} \times_{\mc{Z}} \mc{Z}'$. The given
  $\mc{G}$ action on $\mc{X}$ and the trivial action on $\mc{Z}'$
  induces a free action of $\mc{G}$ on $\tilde{\mc{X}}$. By replacing
  $\mc{X}$ by $\tilde{\mc{X}}$ we see that it suffices to prove that
  $\ed(X) \geq e$.

  If $X \to Y$ is a compression of the $G$-action on $X$ then using
  Lemma \ref{l:descent} we get a smooth scheme $\mc{Y}'$ over $\T$
  with a $\mc{G}$-action such that $Y' \subset \mc{Y}'$ compatibly
  with the $G$-action on $Y$. By Lemma \ref{l:mult}, $\mc{H}_0$ acts
  freely on a non-empty open subset of $\mc{Y}_0$, so the stabilizer
  in $\mc{G}_0$ of a general point of $\mc{Y}_0$ must intersect
  $\mc{H}_0$ trivially. The proposition then follows immediately from
  Lemma \ref{l:stab}.
\end{proof}

\begin{ques} \label{q:ext}%
  Can one characterize all finite $p$-groups $G$ for which there
  exists a finite flat \emph{connected} group scheme $\mc{G}$ over
  some $\T$ such that $G = \mc{G}_K$?
\end{ques}

The technical lemma below will be used (via Lemma \ref{l:sub}) to
check that the dimension hypothesis of Proposition \ref{p:incrit}
holds in the proof of Theorem \ref{t:main}.
\begin{lem} \label{l:quot}%
  Let $k$ be any field and let $G$ be an affine group scheme over
  $k$. Let $H$ be a closed subgroup scheme of $G$ and let
  $\pi:G \to G/H$ be the quotient map. Let $Z$ be an affine scheme
  over $k$ such that $\pi$ factors through
  a   morphism $\pi':G \to Z$. If the scheme-theoretic image $\pi'(H)$ of
  $H$ in $Z$ is a reduced $k$-rational point $z$, then the map on
  tangent spaces $T_z\,Z \to T_{[H]} (G/H)$ is surjective.
\end{lem}

\begin{proof}
  Let $A \sr{f}{\lr} B \sr{g}{\lr} C$ be the maps of $k$-algebras
  corresponding to the maps $G \sr{\pi'}{\lr} Z \to G/H$ and let
  $h = g f$.  Let $m_A$ be the maximal ideal of $C$ corresponding to
  $[H]$, $m_B$ the maximal ideal of $B$ corresponding to $z$ and
  let $I_H \subset C$ be the ideal of $H$. The quotient map $\pi$ is
  faithfully flat and $h(m_A)C = I_H$. We need to prove that the map
  $m_A/m_A^2 \to m_B/m_B^2$ induced by $f$ is injective.

  The assumption on the map $\pi'$ implies that $g(m_B)C \subset
  I_H$. Thus, if an element $x \in m_A$ is such that $f(x) \in m_B^2$,
  then $h(x) \in I_H^2$. Since $h$ is faithfully flat, the map
  $A/m_A^2 \to C\otimes_A (A/m_A^2) = C/I_H^2$ induced by $h$ is
  injective. This implies that $x \in m_A^2$ as desired.
\end{proof}

\begin{rem} \label{r:quot}%
  We have assumed that $G$ is affine only for convenience since this
  is the case we will need later: it is easy to see that the lemma can
  be extended to any $G$ (and $Z$) of finite type over $k$.
\end{rem}

\subsection{Almost multiplicative group schemes} \label{s:ggs}%

As already noted above, in any non-trivial application of Proposition
\ref{p:incrit} the group scheme $\mc{G}$ has to be noncommutative. In
this subsection we discuss a special class of commutative group
schemes for which we can also overcome the non-freeness problem
mentioned in Remark \ref{r:descent}.

\begin{lem} \label{l:split}%
  Let $\mc{G}$ be a finite flat group scheme over $\T$ such that
  $G \cong (\Z/p\Z)^r$. Then $\mc{G} \cong \mc{G}' \times \mc{G}''$,
  where $\mc{G}'$ is connected and $\mc{G}'' \cong (\Z/p\Z)^s$ for some
  $s \leq r$.
\end{lem}

\begin{proof}
  We let $\mc{G}'$ be the connected component of $\mc{G}$ containing
  the image of the identity section $\T \to \mc{G}$. It is clearly a
  closed connected subgroup scheme of $\mc{G}$ which is finite flat
  over $\T$. The group scheme $\mc{G}/\mc{G}'$ is then etale and the
  structure of $G$ implies that it is isomorphic to the constant group
  scheme $(\Z/p\Z)_{\T}^s$ for some $s \leq r$. We let $\mc{G}''$ be
  any closed subgroup scheme of $\mc{G}$ such that the induced map
  $G'' \to (\mc{G}/\mc{G}')_K$ is an isomorphism. Such a group scheme
  exists by the structure of $G$ and since $\mc{G}/\mc{G}'$ is etale
  it follows that the map $\mc{G}'' \to \mc{G}/\mc{G}'$ is an
  isomorphism, so $\mc{G}'' \cong (\Z/p\Z)^s$. The map
  $\mc{G}' \times \mc{G}'' \to \mc{G}$ induced by the two inclusions
  gives the desired isomorphism.
\end{proof}

\begin{lem} \label{l:rank}%
  Let $\mc{G}$ be a finite flat connected group scheme over $\T$
  such that $\mc{G} \cong \mc{G}' \times \mu_p^{r-1}$, where
  $G' \cong \Z/p\Z$.
  \begin{enumerate}
  \item If there exists a map $\mc{G} \to \mc{H}$, with $\mc{H}$ also
    finite flat, which is an isomorphism on generic fibres, then
    $\mc{H} \cong \mc{H}' \times \mu_p^{r-1}$ where
    ${H}' \cong \Z/p\Z$. In particular, $t(\mc{H}) = r$.
  \item If $\mc{G}$ acts faithfully on a smooth scheme $\mc{X}$ with
    $\mc{X}_0$ irreducible, then $\mc{G}_0$ acts generically freely on
    $\mc{X}_0$.
  \end{enumerate}
  
\end{lem}

\begin{proof}
  We have $\mc{G}^D \cong (\mc{G}')^D \times (\Z/p\Z)^{r-1}$ and since
  $\mc{H}^D$ maps generically isomorphically to $\mc{G}^D$, it follows
  by applying Lemma \ref{l:split} to $\mc{H}^D$ and dualising that
  $\mc{H} \cong \mc{H}' \times \mu_p^{r-1}$ where $\mc{H}'$ is
  connected (since we have a map $\mc{G} \to \mc{H}$) and
  ${H}' \cong \Z/p\Z$. Since $\mc{H}'$ is connected and
  $H' \cong \Z/p\Z$, it follows from the classification of connected
  group schemes of order $p$ over $k$ that $t(\mc{H}) = 1$, so
  $t(\mc{H}) \geq 1 + (r-1) = r$. This proves (1).

  For (2), we use that $\mc{G}'_0$ is isomorphic to $\mu_p$ or
  $\alpha_p$ so $\mc{G}_0$ is isomorphic to $\mu_p^r$ or
  $\alpha_p \times \mu_p^{r-1}$. In either case, $\mc{G}_0$ has only
  finitely many subgroup schemes (bounded independently of $k$, as may
  be seen easily using Cartier duality) and (2) follows from this
  since the stabilizer of a general point of $\mc{X}_0$ must be
  constant.

\end{proof}

\begin{prop} \label{p:rank}%
  Let $G = (\Z/p\Z)^r$ and suppose $G$ acts faithfully on a smooth
  quasiprojective variety $X$ over $K$. Let $\mc{G}$ be a finite flat
  model of $G$ over $\T$ and suppose it acts on $\mc{X}$, a
  smooth quasiprojective model of $X$, with $\mc{X}_0$ geometrically
  irreducible.
  \begin{enumerate}
  \item If $\mc{G} \cong \mc{G}' \times \mu_p^{r-1}$ with
    $\mc{G}'$ connected, then $\ed(X) = r$.
  \item If the $\mc{G}$ action is also free then $\ed(X;p) = r$.
  \end{enumerate}
\end{prop}

\begin{proof}
  We may reduce (2) to (1) by Abhyankar's lemma as in the proof of
  Proposition \ref{p:incrit}, so it suffices to prove (2).
  
  Using the definition of essential dimension and Lemma
  \ref{l:descent}, if $\ed(X) =s < r$, we may find a smooth affine
  scheme $\mc{Y}$ over $\T$ with $\mc{Y}_0$ geometrically irreducible
  on which $\mc{G}$ acts with the $G$ action on $Y$ being
  (generically) free.

  We let $\mc{H}$ be the effective model for the action of $\mc{G}$ on
  $\mc{Y}$ given by \cite[Theorem A (ii)]{romagny-eff}; to see that
  this applies we can use \cite[Theorem B]{romagny-eff} since one
  easily checks that the assumptions on $\mc{X}$ imply that it is pure
  over $\T$ in the sense of \cite[Definition 2.1.1]{romagny-eff}. The
  group scheme $\mc{H}$ acts faithfully on $\mc{Y}$ and there is a
  morphism $\mc{G} \to \mc{H}$ which is an isomorphism on generic
  fibres. It then follows from Lemma \ref{l:rank} that $\mc{H}$ acts
  generically freely on $\mc{Y}_0$ and so
  $s = \dim(Y) = \dim(\mc{Y}_0) \geq t(\mc{H}) = r$, a contradiction.
\end{proof}

\begin{rem} \label{r:it}%
  Proposition \ref{p:rank} does not extend in a simple way to more
  general integral models $\mc{G}$ of $(\Z/p\Z)^r$, as shown by
  Example \ref{e:it}. It would be very interesting to find (or even
  classify all) other $\mc{G}$ for which (1) and (2) hold.
\end{rem}

\section{Automorphisms of (truncated)
  \texorpdfstring{$\mu$}{mu}-ordinary $p$-divisible groups}
\label{s:mu}

The goal of this section is to prove some basic facts about
automorphisms of truncated $\mu$-ordinary $p$-divisible groups (with
extra structure) which will be the key to applying the results of
\S\ref{s:gsa} to (certain) unitary Shimura varieties. We assume
throughout that $p>2$.

\subsection{Endomorphisms of \texorpdfstring{$\mu$}{mu}-ordinary
  (truncated) $p$-divisible groups}
\label{s:endord}
Let $\k$ be a finite extension of $\F_p$ and $W(\k)$ the ring of Witt
vectors over $\k$. In this subsection, we recall some basic facts and
results about the building blocks of $\mu$-ordinary group schemes and
$p$-divisible groups from \cite{moonen-serre-tate} and then prove a
lemma (Lemma \ref{l:hom}) about morphisms between such group
schemes. This lemma is the basis for our compution of automorphism
groups of more general $\mu$-ordinary group schemes in the next two
subsections. To describe these objects we use contravariant Dieudonn\'e
theory as in \cite{moonen-serre-tate}.

Let $I$ be the set of homomorphisms from $\kappa$ to $k$ (which we
assume is algebraically closed in this section). This set has a
natural cyclic structure induced by the Frobenius automorphism on $\k$
and we denote the successor of $\tau \in I$ by $\tau+1$. To any
function $\f: I \to \{0,1\}$ we associate a Dieudonn\'e module $M(\f)$
over $W(k)$ with basis $\{e_\tau\}_{\tau \in I}$ with a $W(\k)$-action
defined by $x\cdot e_\tau = \tau(x)e_\tau$ for $x \in W(\k)$. We
define the Frobenius and Verschiebung on $M(\f)$ by defining them on
basis elements as follows:
\begin{equation} \label{e:FV}
  F(e_\tau) =
  \begin{cases}
    e_{\tau+1} & \text{if $\f(\tau) = 0$}, \\
    p\cdot e_{\tau+1} & \text{if $\f(\tau) = 1$},
  \end{cases}
  \ \ \ \ \ \ \
  V(e_{\tau+1}) =
   \begin{cases}
    p\cdot e_\tau & \text{if $\f(\tau) = 0$}, \\
    e_\tau & \text{if $\f(\tau) = 1$}.
  \end{cases}
 \end{equation}
 The Diedonne module $M(\f)$ corresponds to a $p$-divisible group
 over $k$ with an action of $W(\k)$ which we call $\X(\f)$ and we
 denote the $p^n$-torsion of this $p$-divisible group by
 $\X(\f)_n$. By \cite[Corollary 2.1.5]{moonen-serre-tate},
 $\X(\f)$ has a canonical lift to $W(k)$ which we denote by
 $\X^{\can}(\f)$ and we denote its $p^n$-torsion by
 $\X^{\can}(\f)_n$.

 Let $\f^j: I \to \{0,1\}$, $j=1,2$, be two functions, so we have
 $\X^{\can}(\f^j)_n$, $j=1,2$. Our first goal is to compute the
 group scheme over $W(k)$ which represents the functor
 $\Hom(\X^{\can}(\f^1)_n,\X^{\can}(\f^2)_n)$ on
 $W(k)$-algebras given by
 $A \mapsto
 \hom(\X^{\can}(\f^1)_{n/A},\X^{\can}(\f^2)_{n/A})$; here, and
 later, we always assume that all maps are compatible with the
 $W(\k)$-action. Given $\f^j$, $j=1,2$, we define (following
 \cite[2.1.5]{moonen-serre-tate}) $\f':I \to \{0,1\}$ by
 \begin{equation} \label{e:f'}
   \f'(\tau) = \begin{cases}
     0 & \text{if $\f^1(\tau) = \f^2(\tau)$}, \\
     1 & \text{if $\f^1(\tau) \neq \f^2(\tau)$}.
   \end{cases}
 \end{equation}
 We also say that $\f^1 \leq \f^2$ (resp.~ $\f^1 <
 \f^2$) if
 $\f^1(\tau) \leq \f^2(\tau)$ (resp.~$\f^1(\tau) \leq \f^2(\tau)$)
 for all $\tau \in I$ 

 \begin{lem} \label{l:hom} Assume $\f^1 \leq \f^2$ or
   $\f^2 \leq \f^1$. Then
   $\Hom(\X(\f^1)_n, \X(\f^2)_n) = 0$ unless
   $\f^1 \leq \f^2$. If this holds then
   \begin{equation} \label{e:1}%
     \Hom(\X(\f^1)_n,\X(\f^2)_n) = \X(\f')_n .
   \end{equation}
   Furthermore,
   if $\f^1 < \f^2$ then
   \begin{equation} \label{e:2}%
   \Hom(\X^{\can}(\f^1)_n,\X^{\can}(\f^2)_n) =
   \X^{\can}(\f')_n .
   \end{equation}
 \end{lem}
 The lemma is a mild strengthening of \cite[Remarks
 2.3.4(ii)]{moonen-serre-tate} and the proof uses the same methods.
\begin{proof}
  We first note that if $\f^1 \leq \f^2$, then
  $X(\f^1) \times X(\f^2)$ is an ordinary $p$-divisible group with
  $W(\kappa)$-structure in the sense of \cite[\S1]{moonen-serre-tate},
  so we may use all the results proved therein.

  Now note that since $\Hom(\X(\f^1)_n,\X(\f^2)_n)$ is representable
  by a group scheme of finite type over $k$, it is determined by its
  points in Artin local $k$-algebras (with residue field $k$). We use
  the exact sequence (of sheaves in the fppf topology over $\spec(k)$)
  \[
    0 \to \X(\f^1)_n \to \X(\f^1) \sr{p^n}{\lr} \X(\f^1) \to 0
  \]
  which for any $k$-algebra $A$ gives rise to a long exact sequence
  \begin{multline} \label{e:3} 0 \to \Hom(\X(\f^1), \X(\f^2))(A)
    \sr{p^n}{\lr} \Hom(\X(\f^1), \X(\f^2))(A) \to \Hom(\X(\f^1)_n,
    \X(\f^2))(A) \\
    \to \Ext^1(\X(\f^1), \X(\f^2))(A) \sr{p^n}{\lr} \Ext^1(\X(\f^1),
    \X(\f^2))(A) \to \dots.
  \end{multline}
  On the category of Artin local $k$-algebras we have from (the proof
  of) \cite[Theorem 2.3.3]{moonen-serre-tate} that
  $\Ext^1(\X(\f^1), \X(\f^2))$ is the trivial sheaf unless
  $\f^1 < \f^2$ in which case it is represented by $\X(\f')$. Clearly
  $\Hom(\X(\f^1)_n, \X(\f^2))(A) = \Hom(\X(\f^1)_n, \X(\f^2)_n)(A)$
  and by \cite[Corollaire 4.3 a)]{illusie-bt} we see that
  $\Hom(\X(\f^1), \X(\f^2))(A) = 0$ whenever $\f^1 \neq \f^2$.

  If $\f := \f^1 = \f^2$ then $\Hom(\X(\f), \X(\f))(k) = W(k)$ by
  Lemma \ref{l:end}.  We claim that in fact
  $\Hom(\X(\f), \X(\f))(A) = W(k)$ for any Artin local $k$-algebra
  $A$. To see this we note that there is a tautological map
  $W(k) \to \Hom(\X(\f), \X(\f))(A)$ for any such $A$, so it suffices
  to show that any element of $\Hom(\X(\f), \X(\f))(A)$ which reduces
  to $0$ modulo the maximal ideal of $A$ must be $0$, but this again
  follows from \cite[Corollaire 4.3 a)]{illusie-bt}. It then follows
  from \eqref{e:3} that $\Hom(\X(\f)_n, \X(\f)_n)$ is represented by
  the constant group scheme $W_n(\k)$ over $k$. Since this is
  $\X(\mg{0})_n$, where $\mg{0}$ is the zero function on $I$, and
  $\f'$ is also the zero function if $\f^1 = \f^2$, it follows that
  all the statements up to \eqref{e:1} hold.

  We now prove that \eqref{e:2} also holds. We can run through the
  above argument with $\X$ replaced by $\X^{\can}$ and $A$ running
  over Artin local $W(k)$-algebras (all the results we have used from
  \cite{moonen-serre-tate} and \cite{illusie-bt} still hold in this
  context) to conclude that $\X^{\can}(\f')_n$ is a closed subscheme
  of the group scheme representing
  $\Hom(\X^{\can}(\f^1)_n,\X^{\can}(\f^2)_n)$ over $W(k)$. By
  \eqref{e:1} we know that this inclusion induces an isomorphism on
  special fibres. On the other hand, it is clear that the generic
  fibre of $\Hom(\X^{\can}(\f^1)_n,\X^{\can}(\f^2)_n)$ is etale of
  order equal to the order $|\kappa|$ (it is a form of $\kappa$ viewed
  as a constant group scheme over $K$) so in fact the inclusion must
  be an equality. We conclude that \eqref{e:2} holds.
\end{proof}
 
\begin{lem} \label{l:end}%
  For any $\f$ and $n>0$, the tautological map
  $W_n(\k) \to \hom(\X(\f)_n, \X(\f)_n)$ is an isomorphism.
\end{lem}

\begin{proof}
  This is an easy computation using Dieudonn\'e modules: Let $M(\f)_n$
  be the Dieudonn\'e module of $\X(\f)_n$. Since the action of $W_n(\k)$
  on each $e_i$ is induced by different embeddings of $\k$ in $k$, one
  easily sees by induction on $n$ that the $W_n(k)$-linear
  endomorphisms of $M(\f)_n$ consists of diagonal matrices (using the
  basis $\{e_i\}_{i \in I}$ with entries in $W_n(k)$). Now since $F$
  on $M(\f)_n$ is $\sigma$-linear and $V$ is $\sigma^{-1}$ linear, the
  fact that for each $i$, either $F(e_i) = e_{i+1}$ or
  $V(e_{i+1}) = e_i$ easily implies that any such endomorphism must
  come from an element of $W_n(\kappa)$.
\end{proof}

\begin{rem} \label{r:ff}%
  We recall here for later use the standard fact (see, e.g., \cite[I,
  Lemma 1.5]{messing-bt}) that for any ($m$-truncated) $p$-divisible
  group $\Y$ over any base the multiplication by $p$ map from $\Y_n$
  to itself induces a faithfully flat map $\Y_n \to \Y_{n-1}$ for all
  $n \leq m$.
\end{rem}

\subsection{Automorphisms of \texorpdfstring{$\mu$}{mu}-ordinary (truncated) $p$-divisible
  groups} \label{s:autononp}

In this section we recall the classification of $\mu$-ordinary
truncated $p$-divisible groups with $W(\k)$-structure from
\cite[\S1]{moonen-serre-tate} and describe their automorphism groups.

We keep the notation from \S\ref{s:endord}. The classification depends
on an integer $d>0$ and a function $\ff:I \to \{0,1,\dots,d\}$.  Let
$M(\ff)$ be the free $W_n(\k)$-module with basis $\{e_{\tau,j}\}$ with
$\tau \in I$ and $j \in \{1,2,\dots,d\}$. Define $F$ and $V$ on
$M(\ff)$ by
\begin{equation} \label{e:nonp}
  F(e_{\tau,j}) =
  \begin{cases}
    e_{\tau + 1,j} & \text{if $ j \leq d - \ff(\tau)$}, \\
    p \cdot e_{\tau +1} & \text{if  $ j > d - \ff(\tau)$}, 
  \end{cases}
  \ \ \ \ \ \ \
  V(e_{\tau+1, j}) =
   \begin{cases}
    p \cdot e_{\tau, j} & \text{if  $ j \leq d - \ff(\tau)$}, \\
    e_{\tau,j} & \text{if  $ j > d - \ff(\tau)$}.
  \end{cases}
\end{equation}

We define a $W_n(\k)$-action on $M(\ff)$ by
$b \cdot e_{\tau,j} = \tau(b)e_{\tau,j}$ where we denote the map
$W_n(\k) \to W_n(k)$ induced by $\tau$ also by $\tau$.  The Diedonne
module $M(\ff)$ corresponds to a $p$-divisible group over $k$ with an
action of $W(\k)$ which we call $\X(\ff)$ and we denote the
$p^n$-torsion of this $p$-divisible group by $\X(\ff)_n$. By
\cite[Corollary 2.1.5]{moonen-serre-tate}, $\X(\ff)$ has a canonical
lift to $W(k)$ which we denote by $\X^{\can}(\ff)$ and we denote its
$p^n$-torsion by $\X^{\can}(\ff)_n$.

If we denote by $M^j(\ff)$ the $W(k)$-submodule of $M(\ff)$ spanned by
the $\{e_{\tau, j}\}$ for $\tau \in I$ then it is clear from the
definitions that it is a direct summand of $M(\ff)$ as a Dieudonn\'e
module and is a module of the form $M(f^j)$ as defined in
\S\ref{s:endord}. Here the function $f^j: I \to \{0,1\}$ is given by
$f^j(\tau) = 0$ if $j \leq d - \ff(\tau)$ and $f^j(\tau) = 1$
otherwise. Moreover, it is clear that $f^j \leq f^{j'}$ if
$j \leq j'$. Thus $\X(\ff)_n$ is a direct product
$\prod_{j \in \{1,\dots,d\}}\X(f^j)_n$ (with $\X(f^j)$ as defined in
\S\ref{s:endord}).

The product decomposition of $\X(\ff)_n$ and $\X^{\can}(\ff)_n$ as
described above is not intrinsic, but we can make it so by grouping
together the factors which are isomorphic. Let
$\sigma: \{1,\dots,d\} \to \{1,\dots,r\}$ be the unique surjective
non-decreasing function such that the $j,j'$ are in a fibre of
$\sigma$ iff $f^j = f^{j'}$. We then set $\X_i(\ff)$ for
$i \in \{1,\dots,r\}$ to be $\prod_{j \in \sigma^{-1}(i)} \X(f^j)_n$
and similarly for $\X^{\can}(\ff)_n$. So we have
\begin{equation} \label{e:slope0}%
  \X(\ff)_n = \prod_{i \in \{1,\dots,r\}} \X_i(\ff)_n\ \ , \ \ 
  \X^{\can}(\ff)_n = \prod_{i \in \{1,\dots,r\}} \X^{\can}_i(\ff)_n .
\end{equation}

For finite $n$, we would like to have a precise description of the
group scheme $\Aut(\X^{\can}(\ff)_n)$ of the automorphisms of
$\X^{\can}(\ff)_n$ preserving the $W(\k)$-action. However, this is not
finite over $W(k)$ in general and we will restrict ourselves to
describing two closed subgroup schemes: the connected component of the
the identity section $\Aut(\X^{\can}(\ff)_n)^0$ which is finite over
$W(k)$ and the special fibre $\Aut(\X(\ff)_n)$ which is finite over
$k$.

Using the decompositions given in \eqref{e:slope0} we may view
$\Aut(\X^{\can}(\ff)_n)$ as an open subscheme of the space of
$r \times r$ matrices with the $(i,i')$ entry being in
$\Hom(\X^{\can}_i(\ff), \X^{\can}_{i'}(\ff))$, the group operation
then being ``matrix multiplication''. Similar statements hold for
$\Aut(\X(\ff)_n)$.

\begin{lem} \label{l:aut0}%
  $ $
  \begin{enumerate}
  \item For all finite $n$,  $\Aut(\X^{\can}(\ff)_n)^0$ is a finite flat
    group scheme over $W(k)$. Its matricial description is given as
    follows:
    \begin{enumerate}
    \item if $i > i'$ then the $(i,i')$-entry is $0$.
    \item the $(i,i)$-entry is the identity of $\Aut(\X^{\can}_i(\ff)_n)$.
    \item if $i < i'$ then the $(i,i')$-entry is
      $\Hom(\X^{\can}_i(\ff)_n, \X^{\can}_{i'}(\ff)_n)$.
    \end{enumerate}
  \item For all finite $n$, $\Aut(\X(\ff)_n)$ is a finite group
    scheme over $k$.  Its matricial description is given as follows:
    \begin{enumerate}
    \item if $i > i'$ then the $(i,i')$-entry is $0$.
    \item the $(i,i)$-entry is $\Aut(\X_i(\ff)_n)$.
    \item if $i < i'$ then the $(i,i')$-entry is
      $\Hom(\X_i(\ff)_n, \X_{i'}(\ff)_n)$.
    \end{enumerate}
  \end{enumerate}
\end{lem}

\begin{proof}
  This follows immediately from the definition of the decomposition
  $\X^{\can}(\ff)_n$ in \eqref{e:slope0} and Lemma \ref{l:hom}: the
  schemes $\Hom(\X^{\can}(f^j), \X^{\can}(f^{j'}))$ are etale if
  $f^{j'} \geq f^j$ so do not contribute to the identity component of
  $\Aut'(\X^{\can}(\ff)_n)$. Thus, the decreasing filtration of
  $\X^{\can}(\ff)$ given by $\prod_{i \geq s} \X^{\can}_i(\ff)$ is
  preserved by $\Aut(\X^{\can}(\ff)_n)^0$ and it must act as the
  identity on the associated graded. Furthermore, it follows from
  Lemma \ref{l:hom} that all the matrix entries are finite and flat
  over $W(k)$, so (1) is proved.

  Part (2) is proved in essentially the same way, using that
  $\Hom(\X(f^j), \X(f^{j'})) = 0$ if $f^{j} > f^{j'}$.
\end{proof}

\begin{rem} \label{r:hom} The group schemes
  $\Hom(\X^{\can}_i(\ff), \X^{\can}_{i'}(\ff))$ and
  $\Hom(\X_i(\ff), \X_{i'}(\ff))$ occurring in Lemma \ref{l:aut0} can
  be determined precisely 
  by using Lemma \ref{l:hom}. In particular, they are always
  $n$-truncated $p$-divisible groups. The structure of
  $\Aut(\X_i(\ff)_n)$ is also easy to determine: it is the constant
  group scheme over $k$ of invertible
  $|\sigma^{-1}(i)| \times |\sigma^{-1}(i)|$ matrices with
  coefficients in $W_n(\tk)$.
\end{rem}

\begin{cor} \label{c:flat}%
  For all $n> 0$ the maps
  \[
    \Aut(\X^{\can}(\ff)_{n+1})^0 \to \Aut(\X^{\can}(\ff)_n)^0, \ \ \
    \Aut(\X(\ff)_{n+1}) \to \Aut(\X(\ff)_n)
  \]
  induced by restriction are surjective maps of group schemes (i.e.,
  they are surjective maps of sheaves in the fppf topology) and are
  faithfully flat.
\end{cor}

\begin{proof}
  Surjectivity in the case of $\X^{\can}(\ff)$ follows from the
  explicit description of these group schemes given in Lemma
  \ref{l:aut0} together with the surjectivity of the maps
  \[ \Hom(\X^{\can}_i(\ff)_{n+1}, \X_i^{\can}(\ff)_{n+1}) \to
    \Hom(\X^{\can}_i(\ff)_{n}, \X_i^{\can}(\ff)_{n})
  \]
  which follows from Lemma \ref{l:hom}. In the case of $\X(\ff)$ we
  also need that the maps
  \[
    \Aut(\X_i(\ff)_{n+1}) \to \Aut(\X_i(\ff)_n)
  \]
  are surjective which is also an easy consequence of Lemma
  \ref{l:hom}.

  The faithful flatness in the case of $\X^{\can}(\ff)$ can be checked
  fibrewise since both the source and target are flat over $W(k)$; the
  statement for the generic fibre being obvious (since it is etale),
  this reduces us to the statement in the case of $\X(\ff)$. Any
  surjective map of affine group schemes over a field is flat by
  \cite[14.2]{waterhouse} so flatness follows.

\end{proof}

\subsection{Automorphisms of polarised \texorpdfstring{$\mu$}{mu}-ordinary (truncated)
  $p$-divisible groups}

\label{s:autord}
In this section we restrict ourselves to the polarised $\mu$-ordinary
group schemes corresponding to unitary Shimura varieties, since these
are the ones for which we will prove lower bounds for the essential
dimension of congruence covers in \S\ref{s:ed}. We begin by
establishing some notation and recalling the definition and
classification of such group schemes from
\cite[\S3]{moonen-serre-tate}.

As before $\k$ is a finite extension of $\F_p$ and $\tk$ is a
quadratic extension of $\k$. The field $\tk$ has a unique involution
$x \mapsto x^{\ast}$ whose fixed field is $\k$. This involution lifts
to $W(\tk)$ with the ring of $\ast$-invariants being $W(\k)$. Let
$\X_n$ be an $n$-truncated $p$-divisible group over any base. We allow
$n= \infty$ in which case $\X$ is simply a $p$-divisible group. We
denote by $c: \X_n \to \X_n^{DD}$ the canonical double duality
isomorphism. By a \emph{polarisation} or \emph{duality} of $\X_n$ we
shall mean an isomorphism (of group schemes or $p$-divisible groups)
$\lambda:\X_n \to \X_n^D$ such that $\lambda = \lambda^D \circ
c$. Such a duality induces an involution $f \mapsto f^{\dagger}$ on
$\en(\X_n)$. Now suppose we have a $W(\tk)$ structure on $\X_n$, i.e.,
a ring homomorphism $\iota: W(\tk) \to \en(\X_n)$. We impose the
compatibility condition on $\lambda$ and $\iota$ that
$\iota(b^{\ast}) = \iota(b)^{\dagger}$ for all $b \in W(\tk)$ and call
such a triple $(\X_n, \iota, \lambda)$ a \emph{polarized
  ($n$-truncated) $p$-divisible group (or $\mr{BT_n}$) with
  $W(\tk)$-structure}.\footnote{These are the ones of type AU in the
  terminology of \cite{moonen-serre-tate} which are the only ones
  we shall consider here.} Note that $\X_n^D$ acquires a
$W(\tk)$-structure by duality and the compatibility condition means
that $\lambda$ is antilinear as a map of group schemes with
$W(\tk)$-structure.

For our purposes we shall need an explicit classification of those
objects over $k$ which Moonen sometimes calls ordinary in
\cite{moonen-serre-tate}, but we shall always call them $\mu$-ordinary
in order to avoid any confusion. We recall this in terms of their
Dieudonn\'e modules from (\cite[\S3.2.3]{moonen-serre-tate}, Case AU).

Let $\wt{I}$ be the set of homomorphisms from $\tk \to k$. The
classification depends on two parameters, a positive integer $d$ and a
function $\ff: \wt{I} \to \Z_{\geq 0}$ such that
$\ff(\tau) + \ff(\bt) = d$ for all $\tau \in \wt{I}$, where for
$\tau \in \wt{I}$, $\bt := \tau \circ \ast$. As in \S\ref{s:endord}
for $I$, the set $\wt{I}$ has a cyclic ordering and we denote the
successor of $\tau$ for this ordering by $\tau + 1$. Let $M(\ff)$ be
the Dieudonn\'e module defined exactly as in \S\ref{s:autononp} but now
with $\tau \in \wt{I}$ instead of $I$. We continue to use the same
notation $\X(\ff)$, $X(\ff)_n$, $\X^{\can}(\ff)$, \dots, as in
\S\ref{s:autononp} for the associated $p$-divisible group, etc., which
are objects with $W(\tk)$-action.

The polarisation $\lambda$ on the $n$-truncated $p$-divisible group
$\X(\ff)_n$ corresponding to $M(\ff)_n$ can be chosen to be the map
deduced from the pairing $\varphi:M(\ff) \times M(\ff) \to W(k)$ given
on basis elements by $\varphi(e_{\tau,j}, e_{\tau',j'}) = 0 $ unless
$\tau' = \bt$ and $j' = d-j+1$ and
$\varphi(e_{\tau,j}, e_{\bt, d-j+1}) = 1$.  The tuple
$(\X(\ff)_n, \iota, \lambda)$ is then an $n$-truncated polarised
$p$-divisible group with $W(\tk)$-structure as defined above.

The polarisation $\lambda$ lifts canonically to a polarisation on
$\X^{\can}(\ff)_n$ which we still denote by $\lambda$, and then
(together with its natural $W(\tk)$-action) $\X^{\can}(\ff)_n$ also
acquires the structure of a polarised truncated $p$-divisible group
with $W(\tk)$-structure.

We have a product decompostion
$\X(\ff)_n \cong \prod_{j \in \{1,\dots,d\}}\X(f^j)_n$ as in
\S\ref{s:autononp} and the definition of $\lambda$ shows that it is a
direct sum of isomorphisms $\X(f^j) \to \X(f^{d-j+1})^D$ for all $j$.
A similar statement holds for $\X^{\can}(\ff)$.

We also have product decompositions as in \eqref{e:slope0} and then
one sees from the definitions that $\lambda$ induces isomorphisms
$\X_i(\ff)_n \to \X_{r-i+1}(\ff)_n$ and the same with
$\X^{\can}(\ff)$.

We now define $\Aut(\X^{\can}(\ff)_n)$ to be the group scheme of all
structure preserving automorphisms of $\X^{\can}(\ff)_n$, i.e,
automorphisms preserving the $W(\tk)$-action as well as the
polarisation. As in \S\ref{s:autononp} we will restrict ourselves to
describing two closed subgroup schemes: the connected component of the
the identity section $\Aut(\X^{\can}(\ff)_n)^0$ which is finite over
$W(k)$ and the special fibre $\Aut(\X(\ff)_n)$ which is finite over
$k$.

We let $\Aut'(\X^{\can}(\ff)_n)^0$ be the group scheme of
automorphisms preserving only the $W(\tk)$-structure, i.e., we ignore
the polarisation, so it is the group scheme denoted by
$\Aut(\X^{\can}(\ff)_n)^0$ in \S\ref{s:autononp}.
We will describe $\Aut(\X^{\can}(\ff)_n)^0$ as a closed subscheme of
$\Aut'(\X^{\can}(\ff)_n)^0$ and $\Aut(\X(\ff)_n)$ as a closed
subscheme of $\Aut'(\X(\ff)_n)$ using the matricial descriptions from
Lemma \ref{l:aut0}. The subgroup scheme we wish to determine consists
of elements $\alpha$ of $\Aut'(\X^{\can}(\ff)_n)^0$ such that
$\lambda \circ \alpha = (\alpha^D)^{-1} \circ \lambda$ and similarly
for $\Aut'(\X(\ff)_n)$. In terms of the involution $\dagger$ on
$\End(\X^{\can}(\ff))$, we wish to determine all automorphisms
$\alpha$ such that $\alpha^{\dagger} = \alpha^{-1}$.

As noted earlier, $\lambda$ restricts to an isomorphism
$\lambda_j:\X^{\can}(f^j)_n \to \X^{\can}(f^{d-j+1})_n^D$ and so also
an isomorphism
$\lambda_i: \X^{\can}_i(\ff)_n \to
\X^{\can}_{r-i+1}(\ff)^D$. (Similarly, we have an isomorphism
$\lambda_i:\X_i(\ff)_n \to \X_{r-i+1}(\ff)_n^D$.) This implies that the
involution $\dagger$ induces isomorphisms
\[
  \dagger_{i,i'}:  \Hom(\X^{\can}_i(\ff)_n, \X^{\can}_{i'}(\ff))_n \to
  \Hom(\X^{\can}_{r-i'+1}(\ff)_n, \X^{\can}_{r-i+1}(\ff)_n)
\]
such that $\dagger_{r-i'+1, r-i+1} \circ \dagger_{i,i'}$ is the
identity of $\Hom(\X^{\can}_i(\ff)_n, \X^{\can}_{i'}(\ff)_n)$ for all
$i \leq i'$. (Similar statements hold for $\X(\ff)$ instead of
$\X^{\can}(\ff)$.)

If $i+i' = r+1$ then $\dagger_{i,i'}$ is an involution of
$\Hom(\X^{\can}_i(\ff)_n, \X^{\can}_{i'}(\ff)_n)$ and in this case we
let $\Hom(\X^{\can}_i(\ff)_n, \X^{\can}_{i'}(\ff)_n)^+$ be the
invariants of $\dagger_{i,i'}$ and
$\Hom(\X^{\can}_i(\ff)_n, \X^{\can}_{i'}(\ff)_n)^-$ be the
anti-invariants, i.e., subgroup scheme on which $\dagger_{i,i'}$ acts
as $-1$. (We use similar notation with $\X(\ff)$ instead of
$\X^{\can}(\ff)$.)
\begin{lem} \label{l:flat}%
  If $ i \leq i'$ and $i+i' = r+1$, consider the endomorphism of
  $\Hom(\X^{\can}_i(\ff)_n, \X^{\can}_{i'}(\ff)_n)$ given by
  $\mr{id} + \dagger_{i,i'}$. Then the kernel (resp.~image) of
  $\mr{id} + \dagger_{i,i'}$ is finite flat over $W(k)$ and equal to
  $\Hom(\X^{\can}_i(\ff)_n, \X^{\can}_{i'}(\ff)_n)^-$
  (resp.~$\Hom(\X^{\can}_i(\ff)_n, \X^{\can}_{i'}(\ff)_n)^+)$ giving
  rise to a split exact sequence of $n$-truncated $p$-divisible groups
  over $W(k)$ (with $W(\k)$-action):
  \[
    0 \to \Hom(\X^{\can}_i(\ff)_n, \X^{\can}_{i'}(\ff)_n)^- \to
    \Hom(\X^{\can}_i(\ff)_n, \X^{\can}_{i'}(\ff)_n) \to
    \Hom(\X^{\can}_i(\ff)_n, \X^{\can}_{i'}(\ff)_n)^+ \to 0 .
  \]
  The analogous statement with $\X(\ff)$ instead of $\X^{\can}(\ff)$
  also holds.
\end{lem}

\begin{proof}
  Given the definition of $\X^{\can}_i(\ff)_n$, it follows from Lemma
  \ref{l:hom} that $\Hom(\X^{\can}_i(\ff)_n, \X^{\can}_{i'}(\ff)_n)$
  is a finite product of group schemes $\X^{\can}(f')_n$ where
  $f':\wt{I} \to \{0,1\}$ is a function such that $f'(\tau) = f'(\bt)$
  for all $\tau$. Furthermore, the involution $\dagger_{i,i'}$ either
  permutes pairs of factors or preserves one factor on which it acts
  by an involution. For a pair of factors which is permuted it is
  clear that the invariants and anti-invariants are isomorphic to
  $\X^{\can}(f')_n$ and we have a direct sum decomposition.

  For a factor which is preserved, by the definitions in \eqref{e:FV}
  we see that at the level of $\X(f')_n$ this involution is given on
  the Dieudonn\'e module $M(f')_n$ by $e_{\tau} \mapsto e_{\bt}$ for all
  $\tau$. The invariants correpond to the Dieudonn\'e module spanned by
  all $e_{\tau} + e_{\bt}$ and the anti-invariants to the one spanned
  by all $e_{\tau} - e_{\bt}$, both of which correspond to a
  truncated $p$-divisible group over $k$ with a $W(\k)$-action (given
  by the function $f'':I \to \{0,1\}$ defined by extending an element
  of $I$ to $\wt{I}$ and then applying $f'$). The geometric generic
  fibre of $\X^{\can}(f')_n$ is etale, isomorphic to the constant
  group $W_n(\tk)$ (with the involution being the map $\ast$). This
  implies that the ranks of the generic and special fibres are the
  same (in both cases), and so the invariants and anti-invariants of
  the action on such a factor are (finite) flat over $W(k)$ and the
  structure of the special fibre shows that it is an $n$-truncated
  $p$-divisible group.

  The first part of the lemma follows immediately from these
  observations. The statement about $\X(\ff)$ follows from this by
  restriction to the special fibre.
\end{proof}

\begin{prop} \label{p:aut2}%
  The group scheme
  $\mc{G} = \mc{G}(\ff)_n := \Aut(\X^{\can}(\ff)_n)^0$ has a
  decreasing filtration by finite flat closed normal subgroup schemes
  $\mc{G} = \mc{G}(1) \supset \mc{G}(2) \supset \cdots \supset
  \mc{G}(s) \supset \cdots \mc{G}(r) = \{1\}$ such that for
  $1 \leq s < r$ we have
  \begin{equation*}
    \mc{G}(s)/\mc{G}(s+1) \cong  \left (\prod_{t = 1}^{\lfloor
        \tfrac{r-s}{2} \rfloor}
      \Hom(\X_{t}^{\can}(\ff)_n, \X_{t+s}^{\can}(\ff)_n)   \right ) \times \mc{P}
  \end{equation*}
  where
  \[
    \mc{P} = \begin{cases}
      \{1\} & \text{ if $r-s$ is even}, \\
      \Hom(\X^{\can}_{\lceil \tfrac{r-s}{2}\rceil} (\ff)_n,
      \X^{\can}_{\lceil \tfrac{r-s}{2}\rceil + s}(\ff)_n)^- & \text{ if
        $r-s$ is odd}.
    \end{cases}
  \]

  The group scheme $\G = G(\ff)_n := \Aut(\X(\ff)_n)$ has a decreasing
  filtration by finite closed normal subgroup schemes
  $\G = G(0) \supset \G(1) \supset \G(2) \supset \cdots \supset \G(s)
  \supset \cdots \G(r) = \{1\}$ such that for $1 \leq s < r$ we have
  \begin{equation*}
    \G(s)/\G(s+1) \cong 
    \left (\prod_{t = 0}^{-1 +\lfloor
        \tfrac{r-s}{2} \rfloor}
      \Hom(\X_{t+1}(\ff)_n, \X_{t+s+1}(\ff)_n)   \right )
    \times P
  \end{equation*}
  where
  \[
    P = \begin{cases}
      \{1\} & \text{ if $r-s$ is even}, \\
      \Hom(\X_{\lceil \tfrac{r-s}{2}\rceil} (\ff)_n, \X_{\lceil
        \tfrac{r-s}{2}\rceil + s}(\ff)_n)^- & \text{ if $r-s$ is odd}.
    \end{cases}
  \]
  Furthermore,
   \begin{equation*}
     \G(0)/\G(1) \cong 
     \left (\prod_{t = 1}^{\lfloor
         \tfrac{r}{2} \rfloor}
       \Aut(\X_{t}(\ff)_n)   \right )
     \times Q
  \end{equation*}
  where
  \[
    Q = \begin{cases}
      \{1\} & \text{ if $r-s$ is even}, \\
      \UAut(\X_{\lceil \tfrac{r}{2}\rceil} (\ff)_n) &  \text{ if $r$
        is odd}.
    \end{cases}
  \]
Here $\UAut$ denotes the group of unitary automorphisms, i.e., the
subgroup of $\Aut$ consisiting of elements $x$ such that $x^{-1} =
x^{\dagger}$.
\end{prop}

\begin{proof}
  We use the matricial representation from Lemma \ref{l:aut0} and
  write elements of $\Aut(\X^{\can}(\ff)_n)^0$ as matrices
  $\Phi = (\phi_{i,i'})$ and the condition we need to be satisfied is
  $\Phi \cdot \Phi^{\dagger} = \mr{id}$. Multiplying out the matrices
  gives relations:
  \begin{equation} \label{e:phi}
    \sum_{k=1}^r \phi_{i,k} \phi_{r-i'+1, r-k+1}^{\dagger} = 0
  \end{equation}
  if $i < i'$ (where for ease of notation we have dropped the
  subscript from the $\dagger$).  Using that $\phi_{i,i'} = 0$ if
  $i>i'$ and $\phi_{i,i} = \mr{id}$, we can solve these equations
  inductively in $i'-i$. Doing this we see that the $\phi_{i,i'}$ can
  be chosen to be arbitrary elements of
  $\Hom(\X^{\can}_i(\ff)_n, \X^{\can}_{i'}(\ff)_n)$ if $i < r/2$ and
  $i' < r+1 - i$ and then $\phi_{r+i'-1, r+i-1}$ is uniquely
  determined (since the set of equations is preserved by $\dagger$).
  When $i +i' = r+1$ the equation \eqref{e:phi} becomes
  \[
    \phi_{i,i'} + \phi_{i,i'}^{\dagger} + \sum_{k=i+1}^r \phi_{i,k}
    \phi_{i, r-k+1}^{\dagger} = 0 .
  \]
  The sum is invariant under $\dagger$ and since the map
  $\mr{id} + \dagger_{i,i'}$ is faithfully flat onto the
  $\dagger_{i,i'}$-invariants by Lemma \ref{l:flat}, we can solve for
  $\phi_{i,i'}$, the space of solutions being a torsor over the
  $\dagger_{i,i'}$-anti-invariants.

  We define the filtration $\mc{G}(s)$ to be the subgroup scheme with
  $\phi_{i,i'} = 0$ if $0 < i-i' < s$ and then the above description
  of all elements of $\Aut(\X^{\can}(\ff)_n)^0$ makes the claims about
  the filtration in this case clear.

  For the group $G$, we take $G(1)$ to be the identity component.
  This is precisely the special fibre of the group $\mc{G}$ and so the
  claim about $G$ follows immediately from the already proved claim
  for $\mc{G}$ and the structure of the etale quotient, which follows
  easily from Lemma \ref{l:hom}.
\end{proof}

\begin{cor} \label{c:flat1}%
  For all $n> 0$ the maps
  \[
  \Aut(\X^{\can}(\ff)_{n+1})^0 \to \Aut(\X^{\can}(\ff)_n)^0, \ \ \
  \Aut(\X(\ff)_{n+1}) \to \Aut(\X(\ff)_n)
\]
induced by restriction are surjective maps of group schemes (i.e.,
they are surjective maps of sheaves in the fppf topology) and are
faithfully flat.
\end{cor}

\begin{proof}
  Surjectivity in the case of $\X^{\can}(\ff)$ follows by induction
  using the filtration of these group schemes given in Proposition
  \ref{p:aut2} together with the surjectivity of the maps
  \[ \Hom(\X^{\can}_i(\ff)_{n+1}, \X_i^{\can}(\ff)_{n+1}) \to
    \Hom(\X^{\can}_i(\ff)_{n}, \X_i^{\can}(\ff)_{n})
  \]
  which follows from Lemma \ref{l:hom}. In the case of
  $\X(\ff)$ we also need that the maps
  \[
    \Aut(\X_i(\ff)_{n+1}) \to \Aut(\X_i(\ff)_n)
  \]
  are surjective which is also an easy consequence of Lemma
  \ref{l:hom}.

  The faithful flatness can be proved in the same way as in Corollary
  \ref{c:flat}.

\end{proof}

\subsection{The main example} \label{s:ex}%

We now analyze in detail the structure of the automorphism group for a
particular class of examples that will be essential for our
incompressibility result for unitary Shimura varieties (Theorem
\ref{t:main}). We assume that $\k = \F_p$ (so $\tk = \F_{p^2}$) and
$d= 2\d +1$ is odd. $\wt{I}$ consists of two elements which we call
$\tau$ and $\bt$. Define the function
$\ff: \wt{I} \to \{0,1,\dots,d\}$ by $\ff(\tau) = \d$ and
$\ff(\bt) = \d + 1$.

One easily checks that the functions $f^j: \wt{I} \to \{0,1\}$,
$j \in \{1,2,\dots,d\}$, associated to $\ff$ in \S\ref{s:autononp} are
given as follows:
\begin{itemize}
\item $f^j$ is identically $0$ if $j \leq \d$,
\item $f^{\d+1}(\tau) = 0$ and $\f^{\d+1}(\bt) = 1$,
\item $f^j$ is identically $1$ if $j > \d+1$.
\end{itemize}
Thus $\X^{\can}(f^j)_1 \cong \tk$ (viewed as a constant group scheme
with the tautological $\tk$-action) for $j \leq \d$ and
$\X^{\can}(f^j)_1 \cong (\tk)^D$ for $j > \d +1$. To simplify notation
in what follows, we will denote the group scheme
$\X^{\can}(f^{\d+1})_1$ by $\mc{E}$. It follows from the description
of its Dieudonn\'e module that $\mc{E}_0$ is an extension of the group
scheme $\alpha_p$ by itself, and there is a unique subgroup scheme
$\mc{F}_0$ of $\mc{E}_0$ of order $p$ which is the kernel of the
Frobenius map.  In particular, $t(\mc{E}) = 1$. (It is well-known that
$\mc{E}_0$ is isomorphic to the $p$-torsion subscheme of any
supersingular elliptic curve over $k$.)

It follows from the description of the $f^j$ that the integer $r$
associated to $\ff$ defined in \S\ref{s:autononp} is $3$. By
Proposition \ref{p:aut2}, the group scheme
$\mc{G}(\ff)_1 = \Aut(\X^{\can}(\ff)_1)^0$ has a two step filtration,
with a normal subgroup scheme $\mc{H}(\ff)_1$ equal to
$\Hom((\tk)^{\d}, ((\tk)^D)^{\d})^- \subset \Hom((\tk)^{\d},
((\tk)^D)^{\d}) \cong M_{\d \times \d}(\tk^D)$. The involution
$\dagger$ induces an involution on $M_{\d \times \d}(\tk^D)$ which we
also denote by $\dagger$ and is given by applying the non-trivial
automorphism of $\tk$ and taking the anti-transpose, i.e., the
$(i,j)$-entry and the $(\d -j +1, \d -i + 1)$-entry are permuted.  One
sees from this that
$\mc{H}(\ff)_1 \cong (\mu_p)^{\d^2}$.
The quotient $\mc{G}(\ff)_1/\mc{H}(\ff)_1$ is equal to
$\Hom((\tk)^{\d}, \mc{E}) \cong \mc{E}^{\d}$.

\begin{lem} \label{l:sub}%
  Let $\mc{G} = (\mc{G}(\ff)_1)^t \times \mc{H}'$, where $\mc{H}'$ is
  a multiplicative group scheme. Let
  $\mc{H} = (\mc{H}(\ff)_1)^t \times \mc{H}'$, so $\mc{H}$ is a normal
  multiplicative subgroup scheme of $\mc{G}$. Then any subgroup scheme
  $\mc{K}_0$ of $\mc{G}_0$ which intersects $\mc{H}_0$ trivially
  projects to a subgroup scheme of $\mc{G}_0/\mc{H}_0$ which lies in
  the Frobenius kernel. Furthermore, $t(\mc{G}(\ff)_1) = \d^2 + \d$
  and $t(\mc{G}_0/\mc{K}_0) = t(\mc{G}_0) = t(\d^2 + \d) + t(\mc{H}')$.
\end{lem}

\begin{proof}
  The first statement easily reduces to the case that
  $\mc{G} = \mc{G}(\ff)_1$ since Frobenius kernels are compatible with
  products.  We now describe the group scheme $\mc{G}(\ff)_1$ more
  explicitly. By the results of \S \ref{s:autononp} and
  \S\ref{s:autord} it can be represented as a closed subgroup scheme
  of the group scheme given in matricial form as
  \[
    \begin{bmatrix}
      1 & \mc{E}^{\d} & M_{\delta \times \delta}(\tk^D) \\
      0 & 1 & (\mc{E}^D)^{\d} \\
      0 & 0 & 1
    \end{bmatrix}
  \]
  with multiplication being matrix multiplication, using the
  tautological pairing $\mc{E} \times \mc{E}^D \to \tk^D$. We
  drop indices and write the isomorphism
  $\mc{E}^{\d} \to (\mc{E}^D)^{\d}$ induced by the involution
  $\dagger$ simply as $\lambda$. This is a sum of component-wise
  isomorphisms, whose reduction modulo $p$ can be described explicitly
  in terms of Dieudonn\'e modules.

  Using \eqref{e:phi} one sees that the subgroup $\mc{G}(\ff)_1$ is
  then given (in terms of points valued in any $W(k)$-algebra) by
  matrices of the form
  \[
    \begin{bmatrix}
      1 & x & y \\
      0 & 1 & z \\
      0 & 0 & 1
    \end{bmatrix}
  \]
  such that $z = -\lambda(x)$ and
  $y + y^{\dagger} = x\cdot \lambda(x)$.
  We therefore get a scheme theoretic section
  $\sigma$ of the quotient map $\mc{G}(\ff)_1 \to \mc{G}(\ff)_1/\mc{H}(\ff)_1
  \cong \mc{E}^{\d}$ by
  \[
    x \in \mc{E}^{\d} \mapsto \begin{bmatrix}
      1 & x & x \cdot \lambda(x)/2 \\
      0 & 1 & -\lambda(x) \\
      0 & 0 & 1
    \end{bmatrix}
  \]
  since $(x \cdot \lambda(x))^{\dagger} = x \cdot \lambda(x)$ and
  division by $2$ makes sense because $p \neq 2$. This implies that
  the map $T(\mc{G}(\ff)_1) \to T(\mc{G}(\ff)_1/\mc{H}(\ff)_1)$ is
  surjective, so
  $T(\mc{G}(\ff)_1) = T(\mc{H}(\ff)_1) + T(\mc{E}^{\d}) = \d^2 + \d$.
  We also note that for any point $x$ (resp.$~z$) in the Frobenius
  kernel of $\mc{E}_0$ (resp.~of $\mc{E}_0^D$), we have
  $x \cdot z = 1$, the unit element in $\tk^D$, so $\sigma$ restricted
  to $\mc{F}_0^{\d}$ is a morphism of group schemes.
  
  Now suppose that $\mc{K}_0$ is a subgroup scheme of $\mc{G}_0$ which
  interesects $\mc{H}_0$ trivially and consider its image $\mc{K}_0'$
  in $\mc{G}_0/\mc{H}_0$ (which we identify with $\mc{E}_0^{\d}$ using
  the $(1,2)$-component in the matricial description). The inverse map
  $s: \mc{K}_0' \to \mc{K}_0$ can be written as
  \[
    x \in \mc{K}_0' \mapsto  \begin{bmatrix}
      1 & x & x \cdot \lambda(x)/2 + \nu(x)\\
      0 & 1 & -\lambda(x) \\
      0 & 0 & 1
    \end{bmatrix}
  \]
  where $\nu: \mc{K}_0' \to \mc{H}_0$ is a map of $k$-schemes. Letting
  $x,x'$ be any (valued) points of $\mc{K}_0'$, by matrix
  multiplication we see that $s(x)\cdot s(x')$ is equal to
  \[
    \begin{bmatrix}
      1 & x + x' & (x' \cdot \lambda(x')/2 + \nu(x')) - x \cdot \lambda(x') + (x \cdot \lambda(x)/2 + \nu(x))\\
      0 & 1 & -\lambda(x + x') \\
      0 & 0 & 1
    \end{bmatrix}
  \]
  Since $\mc{E}_0^{\d}$, hence also $\mc{K}_0'$, is commutative it
  follows that we must have $x\cdot\lambda(x') = x' \cdot \lambda(x)$
  for all $x,x'$ as above.

  If $\mc{K}_0'$ is not contained in
  $\mc{F}_0^{\d} \subset \mc{E}_0^{\d}$, the structure of $\mc{E}_0$
  shows that it must surject onto one of the factors. Since the
  pairing
  $\mc{E}^{\d} \times (\mc{E}^D)^{\d} \to M_{\delta \times
    \delta}(\tk^D)$ is defined component-wise, the existence of
  $\mc{K}_0$ as above and the commutativity of $\mc{E}_0^d$ would
  imply that for all $x, x' \in \mc{E}_0(\Spec(A))$, where $A$ is any
  $k$-algebra, we have $x \cdot \lambda(x') = x' \cdot \lambda(x)$.
  Here $\lambda_{\mc{E}_0}: \mc{E}_0 \to \mc{E}_0^D$ is the
  isomorphism induced by the polarisation restricted to
  $\mc{E}_0$. However, since $\lambda$ is $\tk$-antilinear, the
  (perfect) pairing $(x,x') \mapsto x \cdot \lambda(x')$ is hermitian
  (with respect to the $\tk$-structure) and not symmetric, so we must
  have $\mc{K}_0' \subset \mc{F}_0^{\d}$ as claimed. We also note that
  any subgroup scheme of $\mc{F}_0^{\d}$ has a unique lift to a
  subgroup of $\mc{G}_0$ (given by the restriction of $\sigma$) since
  $\mc{F}_0^{\d}$ is a unipotent group scheme and $\mc{H}_0$ is
  multiplicative.

  We now prove the last statement about dimensions (so now $\mc{G}$ is
  once again as in the statement of the lemma). By what we have seen
  above, $\mc{K}_0$ projects (isomorphically) into the Frobenius
  kernel of $\mc{Q}_0 :=\mc{G}_0/\mc{H}_0$ and is contained in the
  image of the section $\sigma:\mc{Q}_0 \to \mc{G}_0$ (as defined
  above for $\mc{G}(\ff)_1$ and extended component-wise). Now consider
  the quotient map $q: \mc{G}_0 \to \mc{G}_0/\mc{K}_0$. We have a
  commutative diagram
  \[
  \xymatrix{
    \sigma(\mc{Q}_0) \ar[r] \ar[d] & \mc{X}_0 \ar[d] \ar[dr] &  \\
    \mc{G}_0 \ar[r] &  \mc{G}_0/ \mc{K}_0 \ar[r] & \mc{Q}_0/q(\mc{K}_0)
  }
\]
Here $\mc{X}_0$ is defined to be the scheme-theoretic image of
$\sigma(\mc{Q}_0)$ in $\mc{G}_0/ \mc{K}_0$ and all the other maps are
induced by the inclusions or quotient maps. Now $\mc{K}_0$ is
contained in $\sigma(\mc{Q}_0)$, and since the scheme theoretic image
of $\mc{K}_0$ in $\mc{G}_0/ \mc{K}_0$ is trivial, so is its image in
$\mc{X}_0$. We may thus apply Lemma \ref{l:quot} (with $G$ there being
$\mc{Q}_0$, $H$ there being $\mc{K}_0$, and $X$ there being
$\mc{X}_0$) to conclude that the map on tangent spaces
$T(\mc{X}_0) \to T(\mc{Q}_0/q(\mc{K}_0))$ is surjective. A fortiori,
the map $T(\mc{G}_0/ \mc{K}_0) \to T(\mc{Q}_0/q(\mc{K}_0))$ is also
surjective.

Since $\mc{K}_0 \cap \mc{H}_0$ is trivial, it follows that we have an
injection $T(\mc{H}_0) \to T(\mc{G}_0/ \mc{K}_0)$. As $T (\mc{H}_0)$
maps to $0$ in $T(\mc{Q}_0/q(\mc{K}_0))$, we conclude that
$t(\mc{G}_0/ \mc{K}_0) \geq t(\mc{H}_0) +
t(\mc{Q}_0/q(\mc{K}_0))$. Since $\mc{G}_0/ \mc{K}_0$ is a principal
$\mc{H}_0$ bundle over $\mc{Q}_0/q(\mc{K}_0)$ and
$\mc{H}_0 \cong (\mu_p)^b$ for some $b$, $\mc{G}_0/ \mc{K}_0$ embeds in
a principal $\mb{G}_m^b$-bundle over $\mc{Q}_0/q(\mc{K}_0)$ (which must be
trivial), so it follows that the inequality is actually an
equality. Now $t(\mc{G}_0) = t(\mc{H}_0) + t(\mc{Q}_0)$ (by the
existence of the section $\sigma$) so it suffices to show that
$t(\mc{Q}_0) = t(\mc{Q}_0/q(\mc{K}_0))$. This follows from Lemma \ref{l:e0}
since  $\mc{Q}_0 \cong \mc{E}_0^m$ for some $m$ and $q(\mc{K}_0)$
is contained in its Frobenius kernel.
\end{proof}

\begin{lem} \label{l:e0}%
  Let $m$ be any positive integer and $\mc{E}'_0$ a closed subgroup
  scheme of $\mc{E}_0^m$, contained in $\mc{F}_0^m$, its Frobenius
  kernel. Then $t(\mc{E}_0^m/\mc{E}'_0) = m$.
\end{lem}
\begin{proof}
  Note that the Frobenius map $F$ of $\mc{E}_0^m$ has image
  $\mc{F}_0^m$, so $\mc{E}_0'' := F^{-1}(\mc{E}_0')$ is a closed
  subscheme of $\mc{E}_0^m$ whose image in $\mc{E}_0^m/\mc{E}'_0$ is
  precisely its Frobenius kernel. Since the image has order $p^m$ (as
  the kernel is $\mc{E}_0' \subset \mc{E}_0''$), the lemma follows.
\end{proof}

\section{Shimura varieties}  \label{s:shim}%

\subsection{PEL Shimura varieties} \label{s:pel}%

Let $(\msf{\msf{G}},X)$ be a Shimura datum which we assume to be of
PEL type in the sense of Kottwitz \cite[\S5, \S7]{kottwitz-points}. We
assume that $p \neq 2$, the group $\msf{G}_{\Q_p}$ is unramified and
we choose a compact open subgroup
$C = C_p \times C^p \subset \msf{G}(\Q_p) \times \msf{G}(\A_f^p) =
\msf{G}(\A_f)$ such that $C_p$ is hyperspecial. We assume that $C^p$
is small enough and $\R$ is large enough so that by
\emph{op.~cit.}~the Shimura variety corresponding to this data has a
smooth integral model $\mc{S}_C$ over $\T$ (corresponding to some
embedding of the reflex field $E$ of the Shimura datum into the field
$K$). Henceforth, the compact subgroup $C$ will be fixed and since
none of our statements wil depend on the choice of $C$ (assuming $C_p$
is hyperspecial) we drop it from the notation, denoting $\mc{S}_C$
simply by $\mc{S}$.

The scheme $\mc{S}$ carries a universal family of abelian varieties
$\mc{A} \to \mc{S}$ with (certain extra structures) up to prime-to-$p$
isogeny. Its special fibre $\mc{S}_0$ has a dense open subscheme
$\mc{S}_0^{\ord}$ called the $\mu$-ordinary locus
\cite{wedhorn-ordinariness}; its points can be characterised in terms
of the Newton polygon of the associated $p$-divisible group or even
just the structure of its $p$-torsion (with the extra structure)
\cite{moonen-serre-tate}.  If $k$ is algebraically closed, for all $k$
points in the $\mu$-ordinary locus the associated $p$-divisible group
with extra structure is isomorphic and we denote it by $\X$. For any
$n>0$, we denote the $p^n$-torsion of $\X$ by $\X_n$.

Let $\mr{Ig}^{\ord}$ be the Igusa ``variety'' over $k$ as in
\cite[Definition 4.3.1]{cs-compact-unitary} with the $b$ there
corresponding to the $\mu$-ordinary locus (we drop the $\mu$ for
notational convenience).
It is the inverse limit of the schemes $\mr{Ig}^{\ord}_n$ (so $\mr{Ig}^{\ord}$ is
not really a variety) over $\mc{S}_0^{\ord}$ parametrizing
isomorphisms of $\mc{A}[p^n]|_{\mc{S}_0^{\ord}}$ with
$\X_n\times_k \mc{S}_0$ compatible with all extra structures. The maps
$\mr{Ig}^{\ord}_n \to \mc{S}_0$ are all affine morphisms hence so is the map
$\mr{Ig}^{\ord} \to \mc{S}_0^{\ord}$.

Let $\Aut(\X_n)$ denote the group scheme of automorphisms of $\X_n$
(preserving the extra structure). This is always an affine group
scheme and if $n=1$ it follows from \cite[Theorem 2.1.2
(ii)]{moonen-dimension} that it is a finite group scheme.

We would like to know whether $\mr{Ig}^{\ord}_n$ is smooth for all $n$. The
following lemma gives a simple criterion for this which we shall check
holds in certain cases of interest to us.

\begin{lem} \label{l:reduced}%
  If all the maps $\Aut(\X_{n+1}) \to \Aut(\X_{n})$ induced by
  restricting an automorphism of $\X_{n+1}$ to $\X_n$ are faithfully
  flat, then each $\mr{Ig}^{\ord}_n$ is smooth.
\end{lem}

\begin{proof}
  By \cite[Corollary 4.3.9]{cs-compact-unitary}, the map
  $\mr{Ig}^{\ord} \to \mc{S}_0^{\ord}$ is faithfully flat and since it
  is an affine morphism, it is an fpqc covering. This implies that for
  all $n>0$, $\mr{Ig}^{\ord}_n$ is an $\Aut(\X_n)$-torsor over
  $\mc{S}_0^{\ord}$.  If all the maps
  $\Aut(\X_{n+1}) \to \Aut(\X_{n})$ induced by restriction are
  faithfully flat then so are all the maps
  $\mr{Ig}^{\ord}_{n+1} \to \mr{Ig}^{\ord}_n$. Since tensor products
  commute with direct limits, this implies that all the maps
  $\mr{Ig}^{\ord} \to \mr{Ig}^{\ord}_n$ are also faithfully flat.

  By \cite[Corollary 4.3.5]{cs-compact-unitary}, $\mr{Ig}^{\ord}$ is a
  perfect scheme, i.e., the Frobenius map is an automorphism. In
  particular, $\mr{Ig}^{\ord}$ is reduced. Since any faithfully flat map of
  commutative rings must be injective, it follows that $\mr{Ig}^{\ord}_n$ is
  also reduced. Since $k$ is a perfect field, it follows that
  $\mr{Ig}^{\ord}_n$ is generically smooth. By the Serre--Tate theorem, the
  completion of the local ring of $\mr{Ig}^{\ord}_n$ at any $k$-valued point
  only depends on the $p$-divisible group $\X$ (with extra structure),
  so we deduce that $\mr{Ig}^{\ord}_n$, being generically smooth, is smooth
  at all points.

\end{proof}

Let $C_p' \subset C_p$ be the kernel of the reduction map from
$\mc{G}(\Z_p) \to \mc{G}(\F_p)$ and set $C' = C_p' \times
C^p$. Corresponding to $C'$ we have a finite etale cover $S(p)$ of
$S = \mc{S}_K$ (defined over $K$).

Let $\X^{\can}$ be the canonical lift of $\X$ to $W(k)$ and let
$\X_n^{\can}$ denote its $p^n$-torsion. For any $s \in
\mc{S}_0^{\ord}(k)$, $\X^{\can}$ is isomorphic to the $p$-divisible
group of the canonical lift of $s$ which gives a $\T$-valued
point of $\mc{S}$. We replace $\R$ by a finite extension so that
the generic fibre of $\X_1^{\can}$ is a constant group.

Let $\Aut(\X_1^{\can})^0$ be the connected component of the identity
section of the group scheme $\Aut(\X_1^{\can})$. It follows from the
finiteness of $\Aut(\X_1)$ cited earlier
that $\Aut(\X_1^{\can})^0$ is finite over $W(k)$. We do not know if it
is always flat over $W(k)$, but this seems likely.

\begin{prop} \label{p:action}%
  Assume that the hypotheses of Lemma \ref{l:reduced} are satisfied.
  Let $\mc{S}^{\nu}(p)$ be the normalisation of $\mc{S}$ in
  $S(p)$. Then $\mc{S}^{\nu}(p)_0$ is reduced and smooth at all
  $\mu$-ordinary points (i.e., points lying over
  $\mc{S}_0^{\ord}$). There is a rational action of
  $\Aut(\X_1^{\can})$ on $\mc{S}^{\nu}(p)$ such that the restriction
  of the action to the subgroup scheme $\Aut(\X_1^{\can})^0$ is
  regular and free on an open subset $\mc{S}^o(p)$ of
  $\mc{S}^{\nu}(p)$ which dominates $\mc{S}$ and is smooth (and
  surjective) over $\T$.
\end{prop}

\begin{proof}

  Let $\pi:\mc{S}(p) \to \mc{S}$ be the scheme parametrising
  isomorphisms $\mc{A}[p] \to \X_1^{\can} \times_{\T} \mc{S}$
  (preserving the extra structure). The group scheme
  $\Aut(\X_1^{\can})$ is quasifinite over $\T$ by \cite[Theorem 2.1.2
  (ii)]{moonen-dimension} (since its generic fibre is clearly finite),
  so the morphism $\pi$ is quasifinite. The induced morphism
  $\mc{S}(p)_K \to S$ is finite and etale. By the choice of $R$, it
  follows that $\mc{S}(p)_K$ is isomorphic to $S(p)$ as defined
  earlier, justifying our notation. Furthermore, $\mc{S}(p)_0$ is
  equal to $\mr{Ig}^{\ord}_1$ (by definition), so it is smooth over
  $\mc{S}_0$.

  Let $\mc{S}^o(p)$ be the Zariski closure of $S(p)$ in $\mc{S}(p)$
  (with the reduced induced scheme structure). It is flat over $\T$
  and the closed fibre $\mc{S}^o(p)_0$ of $\mc{S}^o(p)$ dominates
  $\mc{S}_0$: for any point $x \in \mc{S}_0^{\ord}(k)$ the canonical
  lift of the abelian variety associated to $x$, i.e., the lift
  corresponding to the lift $\X^{\can}$ of $\X$ by the Serre--Tate
  theorem \cite[Theorem 1.2.1]{katz-st}, gives rise, using
  \cite[Proposition 2.3.12]{moonen-serre-tate}, to a morphism
  $x^{\can}: \T \to \mc{S}$ lifting the point $x$. If we let $\mc{B}$
  be the abelian scheme over $\T$ corresponding to $x^{\can}$, then
  there exists an isomorphism $\mc{B}[p] \to \X_1^{\can}$ (preserving
  the extra structure) by definition of the canonical lift, so
  $x^{\can}$ can be lifted to a map $\wt{x^{\can}}: \T \to \mc{S}(p)$
  whose image lies in $\mc{S}^o(p)$.  Since
  $\dim(\mr{Ig}^{\ord}_1) = \dim(S(p))$, it follows that
  $\mc{S}^o(p)_0$ is a union of irreducible components of
  $\mr{Ig}^{\ord}_1$.
  Then since $\mr{Ig}^{\ord}_1$ is smooth by Lemma \ref{l:reduced}, it
  follows that $\mc{S}^o(p)$ is an open subscheme of $\mc{S}(p)$.

  Let $\mc{S}^{\nu}(p)$ be the normalisation of $\mc{S}$ in $S(p)$. Since
  $\mc{S}^o(p)$ is normal and quasifinite over $\mc{S}$, by Zariski's
  main theorem there exists an open embedding
  $\mc{S}^o(p) \to \mc{S}^{\nu}(p)$ (extending the identity map on the
  generic fibres).  It follows from the smoothness of $\mc{S}^o(p)_0$
  that $\mc{S}^{\nu}(p)_0$ has at least one component which is
  generically smooth. Since $S(p) \to S$ is Galois (by our assumption on
  $\R$), it follows that all components of $\mc{S}^{\nu}(p)_0$ are smooth
  at all $\mu$-ordinary points. Furthermore, the flatness of
  $\mc{S}^{\nu}(p)$ over $\T$ and the reducedness of $S(p)$ imply that
  $\mc{S}^{\nu}(p)_0$ is reduced.

  The tautological action of $\Aut(\X_1^{\can})$ on $\mc{S}(p)$
  induces a rational action of $\Aut(\X_1^{\can})^0$ on
  $\mc{S}^{\nu}(p)$ which is regular on $\mc{S}^o(p)$ (which dominates
  $\mc{S}$ and is smooth and surjective over $\T$) and is free by
  definition of the scheme $\mc{S}(p)$ and the finiteness of
  $\Aut(\X_1^{\can})^0$.
\end{proof}

\subsection{Unitary Shimura varieties} \label{s:unitary}

In this section prove that the hypotheses of Proposition
\ref{p:action} hold in the cases of unitary PEL Shimura varieties at
unramified primes. By this we mean the ones of type A in the sense of
\cite[\S5]{kottwitz-points} at primes $p$ at which the algebraic group
$\msf{G}$ of \S\ref{s:pel} is unramified and the subgroup
$C_p \subset \msf{G}(\Q_p)$ is hyperspecial. For a summary of the data
needed to define these Shimura varieties the reader may also refer to
\cite[\S4.1-4.3]{moonen-serre-tate}

We continue with the notation of \S\ref{s:pel} (except we now assume
that the PEL data is of type A), so we denote the integral model
simply by $\mc{S}$ and the corresponding $\mu$-ordinary $p$-divisible
group by $\X$ and we would like to verify the asssumptions on
$\Aut(\X_n)$ and $\Aut(\X^{\can}_n)^0$.  As explained in
\cite[\S4.3]{moonen-serre-tate}, because of our assumption that $p$ is
unramified, the category of $p$-divisible groups (with extra
structure) associated to the Shimura data is, by Morita theory,
isomorphic to a product of categories of the sort considered in
\S\ref{s:autononp} and \S\ref{s:autord}. Thus,
\begin{equation} \label{e:dec}%
  \Aut(\X_n) \cong \prod_{s=1}^m \Aut(\X(\ff_s)_n) \times
  \prod_{s'=1}^{m'}\Aut(\X(\ff_{s'}')_n)
\end{equation}
where
\begin{itemize}
\item for each $s \in {1,\dots,m}$, there is
  \begin{itemize}
  \item a finite extension $\k_s$ of $\F_p$ and a positive integer $d_s$;
  \item setting $I_s:= \hom(\k_s, k)$, a function
    $\ff_s: I_s \to \{0,1,2,\dots, d_s\}$;
  \end{itemize}
  and then $\Aut(\X(\ff_s)_n)$ is the group scheme considered in
  \S\ref{s:autononp}, and
\item for each $s' \in {1,\dots,m'}$ there is
  \begin{itemize}
  \item a finite extension $\tk_{s'}$ of $\F_p$ of even degree and a
    positive integer $d_{s'}'$
  \item setting $\wt{I}_{s'}:= \hom(\tk_{s'}, k)$, a function
    $\ff_{s'}': \wt{I}_{s'} \to \{0,1,2,\dots, d_{s'}'\}$ satisfying
    the assumption of \S\ref{s:autord},
  \end{itemize}
  and then $\Aut(\X(\ff_{s'}')_n)$ is the group scheme considered
  in \S\ref{s:autord}.
\end{itemize}
An analogous isomorphism also holds in the case of $\Aut(\X^{\can}_n)^0$.

With the isomorphism \eqref{e:dec} at our disposal we are ready to
prove the following:
\begin{cor} \label{c:unitary}%
  Let $\mc{S}$ be the integral model of a PEL Shimura of type A over
  $\T$ where $p$ is a prime at which the corresponding group $G$ is
  unramified and the group $C_p$ is hyperspecial. Let $S(p)$ be the
  principal level $p$ cover of $S$ with Galois group $G(\F_p)$ and let
  $\mc{S}^{\nu}(p)$ be the normalisation of $\mc{S}$ in $S(p)$. Then
  $\mc{S}^{\nu}(p)$ is smooth at all $\mu$-ordinary points of
  $\mc{S}^{\nu}(p)_0$ and the restriction of the rational action of
  $\Aut(\X_1^{\can})$ on $\mc{S}^{\nu}(p)$ to the finite flat subgroup
  scheme $\Aut(\X_1^{\can})^0$ is regular and free on the open subset
  $\mc{S}^o(p)$ of $\mc{S}^{\nu}(p)$ which dominates $\mc{S}$ and is
  smooth (and surjective) over $\T$.
\end{cor}

\begin{proof}
  It follows from Corollaries \ref{c:flat} and \ref{c:flat1} and the
  isomorphism \eqref{e:dec} that in the type A case
  $\Aut(\X_1^{\can})^0$ is finite and flat and the assumptions of
  Proposition \ref{p:action} are satisified, so the corollary follows
  from that proposition.
\end{proof}

\section{Incompressibility of congruence covers} \label{s:ed}%

Let $F_0$ be a totally real number field, let $F$
be a CM extension of $F$ and let $c$ be the non-trivial element of
$\gal(F/F_0)$.  We consider a PEL Shimura datum $(\msf{G},X)$ where
$\msf{G}$ is a group of unitary similitudes corresponding to a
$d$-dimensional Hermitian form $h$ over $F$ or to a central division
algebra $D$ over $F$ of degree $d^2$ with an involution $\ast$ of the
second kind, i.e, it fixes $F_0$ and acts non-trivially on $F$. For
each embedding $\tau:F \to \C$ the Hermitian form or the involution
$\ast$ gives rise to a non-negative integer $n(\tau)$: for a Hermitian
form this is the dimension of the $+$-part of the Hermitian form over
$\C$ induced by $\tau$ and in the division algebra case we refer to
\cite[\S 1]{kottwitz-simple} for the definition. These integers
satisfy $n(\tau) + n(\bt) = d$ for all $\tau$, where
$\bt := c \circ \tau$. 

Let $L$ be the reflex field of the Shimura datum. The group
$\aut(\C/\Q)$ acts on the set of integer valued functions on
$\hom(F,\C)$ and in the cases we consider the reflex field is the
fixed field of the stabilizer in $\aut(\C/\Q)$ of the function
$\tau \mapsto n(\tau)$ (we call this function the \emph{type} of the
PEL datum). The field $L$ is clearly contained in the Galois closure
of $\tau(F)$ (for any $\tau$), so it is either totally real or a CM
field. We let $L_0$ be the maximal totally real subfield of $L$. In
the cases we are most interested in we will always have that
$L \neq L_0$.

We fix the subgroup $C \subset \msf{G}(\A^f)$ as in \S\ref{s:pel} and
we let $S = S_C$ be the corresponding Shimura variety. Let $C' \subset
\msf{G}(\A^f)$ also be as in \S\ref{s:pel}, giving a Shimura variety
$S(p) = S_{C'}$, which is naturally a Galois cover of $S$. Our main
theorem is the following:

\begin{thm} \label{t:main}%
  Let $p \neq 2$ be a prime which splits completely in $F_0$. If
  $d = 2\d+1$ is odd and for each $\tau \in \hom(F,\C)$,
  $n(\tau) \in \{0,\d, \d + 1, d\}$, then the cover $S(p)/S$ is
  $p$-incompressible. More precisely, for each irreducible component
  $S'$ of $S$ and $S(p)'$ of $S(p)$ lying over $S'$, the map
  $S(p)' \to S'$ is incompressible.
\end{thm}

\begin{proof}
  By the definition of a Shimura variety, the derived group of
  $\msf{G}(\mb{R})$ is not compact, so there exists a $\tau$ such that
  $n(\tau) = \d$. This implies that $L \neq L_0$ as then $n(\bt) = \d
  + 1$ so the function $\tau \mapsto n(\tau)$ is not invariant under
  $c$.

  Let $P$ be a prime of $L$ lying over $p$. Since $p$ splits
  completely in $F_0$, it also splits completely in $L_0$ because
  $L_0$ is contained in the Galois closure of $F_0$. If the residue
  field at $P$ is also $\F_p$ then incompressibility follows from
  \cite{FKW} so we assume that it is $\F_{p^2}$. Let $\mc{S}$ be the
  local model over $R$ as in \ref{s:pel}, where $R$ is a dvr
  containing the localisation of $O_L$ at $P$. Since the scheme
  $\mc{S}$ is a moduli space of abelian schemes of PEL type, by a
  result of Wedhorn \cite{wedhorn-ordinariness} (see also
  \cite{moonen-serre-tate}) the $\mu$-ordinary locus is open and dense
  in $\mc{S}_0$, hence there is a corresponding $\mu$-ordinary
  polarised $p$-divisible group with suitable endomorphism
  structure. Up to Morita equivalence (as in \cite[\S
  3.1.2]{moonen-serre-tate}) this is a product of groups of the type
  we have considered in \S\ref{s:autord}, with factors being
  parametrised by the $c$-orbits on the set of simple factors of the
  $\Q_p$-algebra $F \otimes_{\Q} \Q_p$. Our assumption on $p$ implies
  that the simple factors of this semisimple algebra are either $\Q_p$
  itself or the unramified quadratic extension $E$ of $\Q_p$, and $c$
  acts on this by permuting pairs of factors isomorphic to $\Q_p$ and
  acting by the non-trivial involution on the factors isomorphic to
  $E$. Thus our factors are of type AL and AU in the sense of
  \emph{loc.~cit.}, where the endomorphism ring is just
  $\Z_p \times \Z_p = W(\F_p) \times W(\F_p)$ in the AL case and
  $\mc{O}_E = W(\F_{p^2})$ in the AU case. In the first case, since
  the involution permutes the factors we are reduced to considering
  $p$-divisible groups without any polarisation, i.e., as in
  \S\ref{s:autononp}, and in the second case $p$-divisible groups with
  endomorphisms by $W(\F_{p^2})$ and a polarisation as in
  \S\ref{s:autord}.

  We now use the classifications of such $p$-divisible groups which
  are ordinary, so we need to determine the function $\ff$
  corresponding to each factor. In the case AL it follows from the
  assumption on the function $n(\tau)$, that the function $\ff$, which
  is just a number since $I$ is a singleton, is any element of the set
  $\{0,\d,\d+1,d\}$. In the case AU, the same assumption implies that
  the values of $\ff$ on $\wt{I}$ (which is a two element set) are
  either $\{0,d\}$ or $\{\d, \d+1\}$.

  We shall prove the incompressibility by combining Corollary
  \ref{c:unitary} with Proposition \ref{p:incrit}, so we must analyze
  the groups $\Aut(\X^{\can}(\ff)_1)^0$ in each case, as the group
  $\mc{G}$ acting on the $\mc{S}'$ produced by Corollary
  \ref{c:unitary} is the product of all these groups over all $\ff$
  occurring above.

  In case AL, since two factors are interchanged we may assume that
  $\ff$ is $0$ or $\d$. When it is $0$, $\Aut(\X^{\can}(\ff)_1)^0$ is
  trivial. When it is $\d$, then $r=2$ and $\X_1(\ff)$ is etale and
  $\X_2(\ff)$ is multiplicative. The group $\Aut(\X^{\can}(\ff)_1)^0$
  is then isomorphic to $\mu_p^{\d(\d+1)}$ as a special case of the
  Lemma \ref{l:aut0}.

  In the case AU, the function $\ff$ is the one considered in
  \S\ref{s:ex} so we may use the computations we have made there.  In
  particular, we see that $\mc{G}$ above is exactly of the form
  considered in Lemma \ref{l:sub} with $\mc{H}'$ being the product of
  all the multiplicative groups corresponding to the factors of type
  AL and $t$ the number of factors of type AU.

  The dimension of any unitary Shimura variety as at the beginning of
  this section is given by $\sum_{\{\tau, \bt\}} n(\tau)\cdot n(\bt)$
  since each pair $\{n(\tau), n(\bt)\}$ contributes a factor
  $SU(n(\tau), n(\bt))$ to the adjoint group of $\msf{G}(\mb{R})$. It
  follows from this and Lemma \ref{l:sub} that
  $t(\mc{G}) = \dim(S) = \dim(S(p))$. We are now ready to apply
  Proposition \ref{p:incrit}. Let $\mc{G}$ be as above and let
  $\mc{H}$ be the subgroup $\mc{H}(\ff)_1^t \times \mc{H}'$ of
  $\mc{G}$.  Lemma \ref{l:sub} shows that the hypotheses of
  Proposition \ref{p:incrit} are satisfied for the $\mc{G}$-action on
  $\mc{S}^o(p)$ with $e = \dim(S)$ so the theorem is proved.
  \end{proof}

\begin{ex}
  If $F_0 = \Q$, so $F$ is an imaginary quadratic field, then all
  primes $p$ satisfy the condition in \ref{t:main} but since $L = F$,
  the results of \cite{FKW} only apply for the primes $p$ which split
  in $F$. As a conequence, one can add the case of odd special unitary
  groups, albeit only over $\F_p$, to the list of groups occuring in
  \cite[Corollary 4.3.12]{FKW}. These groups are also inaccesible by
  the methods of \cite{bf-fixed} since the Hermitian symmetric space
  in the corresponding Shimura datum $(\msf{G},X)$ is not of tube
  type.
\end{ex}

\begin{rem} \label{r:other}%
  $ $
  \begin{enumerate}
  \item The reason that we restrict the values of $n(\tau)$ and also
    only consider primes split in $F_0$ is that the analogue of Lemma
    \ref{l:sub} does not hold for any $\ff$ except for the one
    considered there. However, one can derive explicit inequalities
    for $t(\mc{G}_0/\mc{K}_0)$ in many other cases and these can be
    used to give non-trivial lower bounds, not obtainable by other
    known methods, for the essential dimension at $p$ of
    $p$-congruence covers of more general unitary Shimura
    varieties.\footnote{This will appear in the Ph.D.~thesis of the
      second-named author.}
  \item In the case the Shimura variety corresponds to a division
    algebra, it is proper. The results of \cite{FKW} also apply to
    proper Shimura varieties---this is one of the main advantages of
    their method compared to the method of \cite{bf-fixed}---but in
    that (ordinary) setting the incompressibility essentially comes
    about from the action of an elementary abelian $p$-group (a
    subgroup of the full Galois group) on the congruence cover
    $S(p)$. In the non-ordinary setting the use of noncommutative
    group schemes seems to be unavoidable.
  \end{enumerate}
\end{rem}

\section{Abelian varieties} \label{s:abelian}%

In this section we apply Proposition \ref{p:rank} to prove some
results towards the following conjecture due to
P.~Brosnan.\footnote{Personal communication}

\begin{conj} \label{c:abelian}%
  Let $A$ be an abelian variety over a field $L$ of characteristic
  zero. Then for all primes $p$, the multiplication by $p$ map
  $[p]: A \to A$ is $p$-incompressible.
\end{conj}

It is clear that this conjecture implies that the multiplication by
$n$ map $[n]:A \to A$ is incompressible for all $n>1$. If $\dim(A) =1$
then the conjecture is trivially true and if $\dim(A) = 2$
incompressibility---but not $p$-incompressibility---can be proved in
an elementary way using the known structure of the automorphism groups
of $\mb{P}^1$ and elliptic curves.  We also note that if $\dim(A)$ is
arbitrary but $A$ is sufficiently generic then the conjecture can be
proved by using results of Gabber \cite[Appendice]{jlct-exposant} (or
from the results below) and Brosnan's conjecture for $A$ a product of
elliptic curves is closely related to a question of
Colliot-Th\'el\`ene \cite[Question 1]{jlct-exposant}.

\smallskip

We first note a simple reduction:

\begin{lem} \label{l:spec}%
  Suppose the field $L$ has a discrete valuation at which $A$ has good
  reduction and the residue characteristic is also zero. If the
  conjecture for a fixed prime $p$ holds for the special fibre of the
  Neron model of $A$, then it also holds for $A$.
\end{lem}

\begin{proof}
  We let $K$ be the completion of $L$ with respect to the given
  discrete valuation.  Let $\mc{G}$ be the finite etale group scheme
  over $\T$ given by $\mc{A}[p]$ (which we may assume is constant by
  increasing $K$), where $\mc{A}$ is the Neron model of $A$ over
  $\T$. By Lemma \ref{l:descent}, given a compression $f:A \dar B$ of
  the $G$-torsor $A$ we get a $\mc{G}$-equivariant surjective flat
  morphism $\mc{A}' \to \mc{B}'$ where $\mc{A}'$ is open in $\mc{A}$,
  $\mc{A}_0'$ is non-empty, $\mc{B}'$ is affine and smooth over $\T$
  and on the generic fibres the map on function fields is equal to
  that induced by $f$. It suffices to show that the $G$ action on
  $\mc{B}_0'$ is generically free.

  Let $O(\mc{B}')$ be the affine ring of $\mc{B}'$. The group $G$ acts
  on it and $O(\mc{B}')^G$ is finitely generated over $R$ since
  $O(\mc{B}')$ is so (by Hilbert's argument). It is then clearly
  normal and $O(\mc{B}')$ is finite over $O(\mc{B}')^G$ and we set
  $\mc{B}'/G := \spec(O(\mc{B}'))^G$. Since $\mc{B}'$ is smooth over
  $\T$ it follows that the local ring of $\mc{B}'/G$ at the
  generic point of its special fibre is a dvr which is unramfied over
  $\R$. It follows that the map of function fields induced by the map
  $\mc{B}'_0 \to (\mc{B}'/G)_0$ must be Galois with Galois group $G$,
  so the $G$-action on  $\mc{B}_0'$ is generically free.

\end{proof}

\begin{rem}
  The proof of the lemma gives a very general statement comparing
  essential dimension of the generic fibre and special fibre of a
  $G$-torsor over a dvr.
\end{rem}

Lemma \ref{l:spec} easily allows one to reduce Conjecture \ref{c:ab}
to the case when $A$ is the base change to $L$ of an abelian variety
over a number field, but we do not use this in the following:

\begin{thm} \label{t:ab}%
  Let $A$ be an abelian variety over a field $L$ of characteristic
  zero. Suppose $L$ has a discrete valuation of prime residue
  characteristic $p$ at which $A$ has semi-stable reduction and the
  special fibre of the Neron model of $A$ associated to this valuation
  has a closed subgroup scheme of the form $(\mu_p)^n$ for some
  $n < d:=\dim(A)$. Then the essential dimension at $p$ of
  $[p]:A \to A$ over any extension field of $L$ is at least $n+1$.
\end{thm}

\begin{proof}
  We may replace $L$ by an extension $K$ which is complete with
  respect to a discrete valuation extending the one on $L$. Let $\R$
  be the ring of integers of $K$ and $k$ the residue field which we
  may assume is algebraically closed. We may also assume that $K$ is
  large enough so that $A[p]$ (over $K$) is a constant group scheme.
  Let $\mc{A}$ be the Neron model of $A$ over $\T$ and let $\mc{A}[p]$
  be its $p$-torsion subscheme. Since $A$ has semi-stable reduction,
  $\mc{A}[p]$ is flat over $\T$, so its identity component
  $\mc{A}[p]^0$ is finite flat over $\T$. The closed fibre
  $\mc{A}[p]^0_0$ is a finite group scheme of order at least $p^{n+1}$
  (since $n < d$).

  By applying Lemma \ref{l:split} to the Cartier dual of
  $\mc{A}[p]^0_0$ (and then dualising again), we see that
  $\mc{A}[p]^0_0$ contains a finite flat subgroup scheme $\mc{G}$ of
  rank $n+1$ such that $\mc{G}_0$ contains a subgroup scheme
  isomorphic to $(\mu_p)^n$. The Neron model $\mc{A}$ is smooth over
  $\T$ and $\mc{A}[p]^0_0$, hence also $\mc{G}$, acts freely on it by
  translation. We may therefore apply Proposition \ref{p:rank} to
  conclude the proof.
\end{proof}

\begin{cor} \label{c:ab}%
  Let $A$ be an abelian variety of dimension $d$ over a field $L$ of
  characteristic zero. Suppose $L$ has a discrete valuation of prime
  residue characteristic $p$ at which $A$ has good reduction and the
  special fibre of the Neron model $\mc{A}$ of $A$ associated to this
  valuation has $p$-rank at least $d-1$.  Then the essential dimension
  at $p$ of $[p]:A \to A$ over any extension field of $L$ is $d$.
\end{cor}

\begin{proof}
  This is an almost immediate consequence of Theorem \ref{t:ab}. If
  the $p$-rank of the special fibre is at least $d-1$, then since the
  special fibre is an abelian variety it must contain a subgroup
  scheme isomorphic to $(\mu_p)^{d-1}$.
\end{proof}

\begin{lem} \label{l:prank}%
  Let $A$ be an abelian variety of dimension $d> 1$ over a number
  field $L$.  The $p$-rank of the reduction of a $g$-dimensional
  abelian variety $A$ at a prime $P$ of $O_L$ above $p$ is least two
  for a set of rational primes $p$ of positive density (depending on
  $A$).
\end{lem}

\begin{proof}
  This is well-known if $d=2$ \cite[VI, Corollary 2.9]{dmos} and a
  similar argument works in general. We first note that the set of
  primes $P$ of $O_L$ lying above rational primes $p$ which split
  completely in $O_L$ is of density one.
  We consider the
  characteristic polynomial 
  \[
    \Phi_P(x) = x^{2d} + a_1(P)x^{2d-1} + a_2(P)x^{2d-2} + \dots + a_{2d}(P)
  \]
  of geometric Frobenius for such primes $P$ acting on the $l$-adic
  Tate module of $A$ for some fixed rational prime $l$. The
  $a_i(P) \in \Z$ and by Weil's theorem $|a_i(P)| \leq c_ip^{i/2}$,
  where $c_i$ is a constant depending only on $d$. If the $p$-rank of
  the reduction of $A$ at $P$ is at most one, then $a_2$ is divisible
  by $p$ so $a_2/p$ lies in a finite set of integers. On the other
  hand, since $a_{2d}(P) = p^d$ for all $P$ as above, we see that
  $a_2(P)^d/a_{2d}(P)$ takes on only finitely many values for such
  $P$.
  
  By enlarging $L$ if necessary, we may assume that the Zariski
  closure $G$ of the image of $\mr{Gal}(\ov{L}/L)$ (acting on the Tate
  module $\otimes \Q_l$) is connected.  Each $a_i$, viewed as a
  function on $\mr{Gal}(\ov{L}/L)$ is the restriction of an algebraic
  function on $G$ to the image of $\mr{Gal}(\ov{L}/L)$, hence so is
  $(a_2)^d/a_d$.  The image of $\mr{Gal}(\ov{L}/L)$ is an open
  subgroup of $G(\Q_l)$ and $G$ contains the homotheties by
  \cite{bogomolov-alg} so $(a_2)^d/a_d$ is not a constant function on
  $G$. Thus, the level sets of the function have measure $0$ in
  $G(\Q_l)$, hence by the Chebotarev density theorem the set of all
  $P$ as above with the reduction of $A$ at $P$ being of $p$-rank at
  most one has density $0$.
\end{proof}

\begin{cor} \label{c:three}%
  Brosnan's conjecture in the case $\dim(A) \leq 3$ holds for a set of
  primes $p$ (depending on $A$) of positive density.
\end{cor}

\begin{proof}
  This follows by combining Lemma \ref{l:prank}, Lemma \ref{l:spec}
  and Corollary \ref{c:ab}.
\end{proof}

\begin{rem} \label{r:2}%
  $ $
  \begin{enumerate}
  \item If $\dim(A) = 2$ then, as already noted earlier,
    incompressibility of $[p]$ holds for all $p$ but we do not know
    whether $[p]$ is $p$-incompressible for all $p$ (but see
    \cite[Proposition 11]{jlct-exposant} for a related result in a
    special case).
  \item The difficulty in proving Brosnan's conjecture in general by
    our method lies in the problem of the (possible) non-freeness of
    the descended $\mc{G}$-action mentioned in Remark
    \ref{r:descent}. We do not actually have such an example when
    $\mc{G}$ is the $p$-torsion of an abelian scheme over $\T$ (or
    even a $1$-truncated $p$-divisible group) and it would be very
    interesting to know whether in this case the conclusion of Lemma
    \ref{l:descent} can be made stronger, i.e., if the $\mc{G}$-action
    is free at a general point of $\mc{X}_0$, is it also free at a
    general point of $\mc{Y}'_0$?
  \item It is natural to extend Brosnan's conjecture to abelian
    varieties over arbitrary fields---it suffices to assume
    algebraically closed---of characteristic $l \neq p$, but our
    methods do not apply: if $\dim(A) = 2$ then incompressibility of
    $[p]$ does hold for all $p \neq l$, but for a fixed $A$ we do not
    know whether $p$-incompressibility holds for even a single $p$ if
    $l > 0$. It would also be interesting to
    consider the case $p = l$, when $[p]:A \to A$ is a torsor under a
    nonreduced group scheme.
  \end{enumerate}
\end{rem}

\bibliographystyle{siam} 

\bibliography{../sources}

\end{document}